\documentclass[12pt]{amsart}
\usepackage{amssymb,amsmath}
\makeatletter
\def\@cite#1#2{{\m@th\upshape\bfseries%
[{#1\if@tempswa{\m@th\upshape\mdseries, #2}\fi}]}}
\makeatother
\theoremstyle{plain}
\newtheorem{thm}{Theorem}[section]
\newtheorem{lem}[thm]{Lemma}
\newtheorem{cor}[thm]{Corollary}

\theoremstyle{definition}
\newtheorem{rem}[thm]{Remark}

\newtheorem{eg}[thm]{Example}

\newcommand{\Prf}{\noindent\textbf{Proof.\ }}
\newcommand{\bx}{\strut\hfill$\blacksquare$\medbreak}
\newcommand{\upbx}{\vspace{-2.5\baselineskip}\newline\hbox{}\hfill$\blacksquare$\newline\
medbreak}

\newcommand{\BH}{{\B(\H)}}

\newcommand{\lip}{\langle}
\newcommand{\rip}{\rangle}
\newcommand{\ol}{\overline}
\newenvironment{sbmatrix}{\left[\begin{smallmatrix}}{\end{smallmatrix}\right]}

\newcommand{\td}{\widetilde}
\newcommand{\upplus}[1][2.3]{\!\rule{0pt}{#1ex}^\oplus}

\DeclareMathOperator*{\wotlim}{\textsc{wot}--lim}

\DeclareMathOperator*{\sotlim}{\textsc{sot}--lim}

\newcommand{\wot}{\textsc{wot}}

\newcommand{\bbC}{{\mathbb{C}}}

\newcommand{\bbN}{{\mathbb{N}}}

 \newcommand{\A}{{\mathcal{A}}}
 \newcommand{\B}{{\mathcal{B}}}

 \newcommand{\E}{{\mathcal{E}}}
 \newcommand{\F}{{\mathcal{F}}}
 
\renewcommand{\H}{{\mathcal{H}}}
 
 \newcommand{\J}{{\mathcal{J}}}
 \newcommand{\K}{{\mathcal{K}}}
\renewcommand{\L}{{\mathcal{L}}}
 \newcommand{\M}{{\mathcal{M}}}
 \newcommand{\N}{{\mathcal{N}}}
\renewcommand{\O}{{\mathcal{O}}}

\renewcommand{\S}{{\mathcal{S}}}
 
 \newcommand{\U}{{\mathcal{U}}}
 \newcommand{\V}{{\mathcal{V}}}
 \newcommand{\W}{{\mathcal{W}}}
 \newcommand{\X}{{\mathcal{X}}}
 \newcommand{\Y}{{\mathcal{Y}}}
 
\newcommand{\eps}{\varepsilon}
\renewcommand{\phi}{\varphi}
\newcommand{\upchi}{{\raise.35ex\hbox{$\chi$}}}
\newcommand{\fA}{{\mathfrak{A}}}
\newcommand{\fB}{{\mathfrak{B}}}

\newcommand{\fI}{{\mathfrak{I}}}
\newcommand{\fJ}{{\mathfrak{J}}}
\newcommand{\fK}{{\mathfrak{K}}}
\newcommand{\fL}{{\mathfrak{L}}}
\newcommand{\fM}{{\mathfrak{M}}}

\newcommand{\fP}{{\mathfrak{P}}}

\newcommand{\fR}{{\mathfrak{R}}}
\newcommand{\fS}{{\mathfrak{S}}}

\newcommand{\fW}{{\mathfrak{W}}}

\renewcommand{\And}{\text{ and }}

\newcommand{\qand}{\quad\text{and}\quad}
\newcommand{\qqand}{\qquad\text{and}\qquad}

\newcommand{\qfor}{\quad\text{for}\quad}
\newcommand{\qforal}{\quad\text{for all}\quad}

\newcommand{\Alg}{\operatorname{Alg}}

\newcommand{\dist}{\operatorname{dist}}
\newcommand{\id}{\operatorname{id}}

\newcommand{\Lat}{\operatorname{Lat}}
\newcommand{\prank}{\operatorname{pure\:rank}}
\newcommand{\ran}{\operatorname{Ran}}
\newcommand{\rank}{\operatorname{rank}}
\newcommand{\re}{\operatorname{Re}}

\newcommand{\spn}{\operatorname{span}}
\newcommand{\Tr}{\operatorname{Tr}}
\begin{document}

\title[Isometric Dilations]%
{Isometric Dilations of non-commuting finite
rank n-tuples}
\thanks{draft, March 10, 1999}
\author[K.R.Davidson]{Kenneth~R.~Davidson}
\thanks{The first author was partially supported by an NSERC grant.}
%
\author[D.W.Kribs]{David~W.~Kribs}
\thanks{The second author was partially supported by NSERC and
OGS scholarships.}
%
\author[M.E.Shpigel]{Miron~E.~Shpigel}
%
\date{}
\begin{abstract}
A contractive $n$-tuple $A=(A_1,\dots,A_n)$ has a minimal joint
isometric dilation $S=(S_1,\dots,S_n)$ where the $S_i$'s are
isometries with pairwise orthogonal ranges.  This determines a
representation of the Cuntz-Toeplitz algebra.  When $A$ acts on a
finite dimensional space, the \wot-closed nonself-adjoint algebra
$\fS$ generated by $S$ is completely described in terms of the
properties of $A$.  This provides complete unitary invariants for the
corresponding representations.  In addition, we show that the algebra
$\fS$ is always hyper-reflexive.  In the last section, we describe
similarity invariants.  In particular, an $n$-tuple $B$ of $d\times d$
matrices is similar to an irreducible $n$-tuple $A$ if and only if
a certain finite set of polynomials vanish on $B$.
\end{abstract}
\maketitle

In \cite{DP1,DP2}, the first author and David Pitts studied a
class of algebras coined {\it free semigroup algebras}.  These
are the \wot--closed (nonself-adjoint) unital operator algebras
generated by an $n$-tuple of isometries with pairwise orthogonal
ranges.  When these ranges span the whole space, the associated
norm-closed self-adjoint algebra is a representation of the Cuntz
algebra.  This nonself-adjoint algebra can contain detailed
information about fine unitary invariants of the corresponding
C*-algebra representation.  Indeed in \cite{DP1} the set of atomic
representations of the Cuntz algebra is completely classified.  On
the other hand, when the ranges span a proper subspace, the
representation contains a multiple of the left regular representation
of the free semigroup on $n$ letters.  The \wot-closed algebra of the
left regular representation is called the {\it non-commutative
analytic Toeplitz algebra}.  This nomenclature is justified by a good
analogue of Beurling's Theorem \cite{Pop_beur,AP,DP1},
hyper-reflexivity \cite{DP1} and the relationship \cite{DP2} between
its automorphism group and the group of conformal automorphisms of
the ball in $\bbC^n$.

The connection with dilation theory derives from a theorem of
Frahzo, Bunce and Popescu \cite{Fra1,Bun,Pop_diln}.  If
$A= (A_1,\dots,A_n)$ is an $n$-tuple of
operators such that $AA^* = \sum_{i=1}^n A_iA_i^* \le I$,
then there is a unique minimal isometric dilation to isometries
$S_i$ on a larger space with pairwise orthogonal ranges.
Popescu \cite{Pop_diln} establishes the analogue of Wold's
decomposition which splits this into a direct sum of a multiple
of the left regular representation and a representation of the
Cuntz algebra.  Moreover, Popescu \cite{Pop_vN} obtains the
non-commutative analogue of von Neumann's inequality in this context.
We mention in passing that there has been recent interest in dilating
commuting $n$-tuples as well \cite{Pop_comm,Arv_III,Arv_euler}.

On the other hand, representations of the Cuntz algebra correspond to
endomorphisms of $\B(\H)$ \cite{Pow,Laca,BJP,BJendII}.  This has
created new interest in classifying these representations up to
unitary equivalence.  The well-known theorem of Glimm \cite{Glimm}
shows that this classification is non-smooth because $\O_n$ is
anti-liminal (or NGCR).  Nevertheless, interesting classes of
representations do lend themselves to a complete analysis.
In \cite{BJit}, Bratteli and Jorgensen introduced a class of
representations which turned out to be a special case of the
atomic representations classified in \cite{DP1} using nonself-adjoint
techniques.  In \cite{BJendII} they introduce a different class
associated to finitely correlated states.  The reader will see a lot
of parallels between their results and ours, though the approach is
quite different.  In the end, they specialize to the subclass of
diagonalizable shifts in order to obtain a classification theorem.  In
this paper, we obtain good unitary invariants for the class of all of
these finitely correlated representations.

The goal of this paper is two-fold.  First we wish to understand
the structure of the free semigroup algebra generated by the
dilation of an $n$-tuple $A$ in terms of information obtained
from the $n$-tuple itself (and the algebra it generates).
In particular, we seek unitary invariants for the associated
C*-algebra representation.  Secondly, we wish to determine whether
these algebras are reflexive and even hyper-reflexive.  In this paper,
we focus on the case in which the $n$-tuple $A$ acts on a finite
dimensional space.  Here we obtain a complete description of
the algebra.  This enables us to decompose the associated
representation as a direct sum of irreducible representations and
obtain complete unitary invariants.  These algebras all turn out to be
hyper-reflexive.

In the last section, we discuss similarity invariants.  One of the
surprising consequences is a complete invariant for an irreducible
$n$-tuple of $d \times d$ matrices up to similarity.  An algorithm for
determining if two $n$-tuples of matrices are similar is provided by
Friedland \cite{Fried}.  But this method rapidly gets complicated.  So
it is not clear whether it is superior to ours.  We find that there
is a finite set of no more than $1+(n-1)d^2$ polynomials $p_j$ so that
another $n$-tuple $B$ is similar to $A$ if and only if $p_j(B)=0$ for
all $j$.  These polynomials are obtained from a computable set of
generators of an ideal of the left regular free semigroup algebra as a
right ideal, which amounts to computing an orthonormal basis for a
certain subspace.  In practice, one only needs generators as a
two-sided ideal, and hence the actual number needed is normally
smaller.

 \section{Background}\label{S:background}

Let $\F_n$ denote the unital free semigroup on $n$ letters
$\{ 1,2,\dots,n \}$, and let $\K_n = \ell^2(\F_n)$ denote the Hilbert
space with basis $\{ \xi_w : w \in \F_n \}$, which is known as
$n$-variable Fock space.  The left regular representation $\lambda$
of $\F_n$ is given by  $\lambda(v) \xi_w := L_v \xi_w = \xi_{vw}$.
In particular, the generators of $\F_n$ determine isometries $L_i$ for
$1 \le i \le n$ with orthogonal ranges such that
$\sum_{i=1}^n L_iL_i^* = I - P_e$
where $P_e = \xi_e \xi_e^*$ is the rank one projection onto the
basis vector for the empty word $e$, which is the identity of
$\F_n$.  The algebra $\fL_n$ is the \wot-closed algebra
generated by the $n$-tuple $L = (L_1,\dots,L_n)$. See
\cite{DP1,DP2,DP3,Kri,Pop_fact,Pop_beur} for detailed information about
this algebra.

More generally if $S_i$, $1 \le i \le n$, are isometries with
$\sum_{i=1}^n S_i S_i^* \le I$,
let $\fS$ denote the unital \wot-closed (nonself-adjoint)
algebra generated by them.  We denote by $S_v$ the isometry
$v(S):= v(S_1,\dots,S_n)$ for each $v \in \F_n$.  A subspace $\W$
is called {\it wandering} for the $n$-tuple $S = (S_1,\dots,S_n)$
provided that the subspaces $S_v \W$ are pairwise orthogonal for
all $v \in \F_n$.  Thus the smallest $\fS$-invar\-iant subspace
containing a wandering space $\W$ is
$\fS[\W] = \sum^\oplus_{v \in \F_n} S_v \W$. The restriction of $\fS$
to this subspace is evidently a multiple of the left regular
representation algebra $\fL_n$, where the multiplicity is given by
$\dim \W$.  Popescu's Wold decomposition \cite{Pop_diln} works as
follows: the subspace $\W = \ran(I - \sum_{i=1}^n S_i S_i^* )$ is
easily seen to be wandering.  Moreover the complement
$\N = \fS[\W]^\perp$ is also invariant for $\fS$, and the restriction
to $\N$ yields isometries $T_i = S_i|_\N$ satisfying
$\sum_{i=1}^n T_iT_i^* = I_\N$.

Suppose $A = (A_1,\dots,A_n)$ is an $n$-tuple of
operators on a Hilbert space $\V$ such that
$A A^* = \sum_{i=1}^n A_i A_i^* \le I$.
Frahzo \cite{Fra1} (for $n=2$), Bunce \cite{Bun} (for $n<\infty$)
and Popescu \cite{Pop_diln} (for $n=\infty$) show that there is
a joint dilation of the $A_i$ to isometries $S_i$ on a
Hilbert space $\H = \V \oplus \K$ which have pairwise orthogonal
ranges.  Popescu observes that if this dilation is minimal in
the sense that $\H = \spn \{ S_v \V : v \in \F_n \}$, then the
dilation is unique (up to a unitary equivalence which fixes
$\V$).  We will always work with this minimal isometric dilation.

Popescu also observes \cite{Pop_vN} that the norm-closed nonself-adjoint
algebra $\A_n$ spanned by $\{L_w : w \in\F_n\}$ is the appropriate
non-commutative analogue of the disk algebra for a version of von
Neumann's inequality.  Namely, if $A$ is a contractive $n$-tuple as
above, then $\|p(A)\| \le \|p(L)\|$ for every non-commuting polynomial
in $n$ variables.  This is immediate from the dilation theorem and the
fact that there is a contractive homomorphism of $\E_n$ onto $\O_n$,
the two possible C*-algebras for the dilation.  However, it turns out
that this quotient map is completely isometric on $\A_n$.  So this norm
estimate is an equality for any contractive $n$-tuple of
{\it isometries}.  This shows that $\O_n$ is the C*-envelope of
$\A_n$.

This presents a rather precise picture for the {\it norm-closed} algebra
generated by an $n$-tuple of isometries with orthogonal ranges.
However, the \wot-closed algebras can be quite different.  They can
reflect the fine unitary invariants of the representation.  The
case $n=1$ is familiar, where the \wot-closed algebra depends on the
spectral invariants of the unitary part and the multiplicity of the
shift (from the Wold decomposition).

When $\sum_{i=1}^n S_iS_i^* = I$, the C*-algebra generated by the
isometries $S_i$ is the Cuntz algebra $\O_n$;  and when
$\sum_{i=1}^n S_iS_i^* < I$,  this C*-algebra is $*$-isomorphic to the
Cuntz--Toeplitz algebra $\E_n$ generated by the left regular
representation $\lambda$.  This algebra is an extension of
the compact operators $\fK$ by $\O_n$.
We associate to each $n$-tuple $S_i$ the
representation $\sigma$ of $\E_n$ given by $\sigma(s_i)=S_i$,
where $s_i$ are the canonical generators of $\E_n$.  When
$\sum_{i=1}^n S_iS_i^* = I$, we may consider this as a representation
of $\O_n$ instead.  Let $\fS_\sigma$ denote the \wot-closed
non-self-adjoint algebra determined by the representation $\sigma$.
One can view the Wold decomposition as the spatial view of the
C*-algebra fact that every representation $\sigma$ of $\E_n$ splits as
a direct sum $\sigma = \lambda^{(\alpha)} \oplus  \tau$ of a
representation $\lambda^{(\alpha)}$, which is faithful on
$\fK$ and thus is a multiple of the identity representation $\lambda$,
and a representation $\tau$ which factors through $\O_n$.

A representation is called {\it atomic} if there is an
orthonormal basis $\{\xi_j\}$ which are permuted up to scalars by
the generating isometries $S_i$.  That is, for each $i$ there is
an endomorphism $\pi_i:\bbN \to \bbN$ and scalars $\lambda_{i,j}$
of modulus 1 such that $S_i \xi_j = \lambda_{i,j} \xi_{\pi_i(j)}$.
These representations decompose as a direct integral of
irreducible atomic representations \cite{DP1}, and these irreducible
atomic representations are of three types.  The first is just the left
regular representation; which is the only one which does not factor
through $\O_n$.  The second type is a class of inductive limits of the
left regular representation, and are classified by an infinite word
(up to shift-tail equivalence) that describes the imbeddings.  The
third type fits into the context of this paper, and so we describe it
in more detail.  See \cite{DP1} for a complete description.

The third type is given by a word $u = i_1 i_2 \dots i_d$ in $\F_n$
and a scalar $\lambda$ of modulus 1.  A finite dimensional space $\V$
of dimension $d$ is formed with a basis $e_1,\dots,e_d$.  Operators
$A_j$, $1 \le j \le n$, are partial isometries given by
\begin{align*}
 A_j e_k &= \delta_{ji_k} e_{k+1} \qfor 1 \le k < d \\
 A_j e_d &= \lambda \delta_{ji_d} e_1 .
\end{align*}
The minimal isometric dilation of this $n$-tuple yields isometries
$S_j$ acting on a space $\H = \V \oplus \K$.
The isometry $S_{i_k}$ maps $e_k$ to $e_{k+1}$ (or $\lambda
e_1$ when $k=d$) and the other $n-1$ isometries send $e_k$ to pairwise
orthogonal vectors which are all wandering vectors for $\fS$.  Thus
$\K = \V^\perp$ is determined by a wandering space $\W$ of dimension
$d(n-1)$, and therefore $\K = \fS[\W] \simeq \K_n^{(d(n-1))}$.
The associated representation $\sigma_{u,\lambda}$ is
irreducible precisely when the word $u$ is primitive, meaning that it
is not a power of a smaller word.  In this case, $\fS$ can be
completely described as the sum of $\B(\H)P_\V$ and a multiple of
$\fL_n$ acting on $\K$ via its identification with $\K_n^{(d(n-1))}$.
The invariant subspaces of this algebra are readily described, and it
turns out to be hyper-reflexive.  (See below).

For the future, we wish to name the type of algebra which occurs here.
Let $\fB_{n,d}$ denote the \wot-closed algebra on a Hilbert space
$\H = \V\oplus\K_n^{(d(n-1))}$ where $\dim \V = d$ given by
\[
  \fB_{n,d} = \B(\H) P_\V + \bigl( 0_\V \oplus \fL_n^{(d(n-1))} \bigr).
\]

Another class of representations which have been studied are the
{\it finitely correlated} representations \cite{BJendII}.  A
representation of $\O_n$ is finitely correlated if there is a finite
dimensional cyclic subspace $\V$ which is invariant for each $S_i^*$.
Likewise, a {\it finite correlated state} is a state $\phi$ such that
in the GNS construction, the invariant subspace for the $S_i^*$'s
generated by the cyclic vector $\xi_\phi$ is finite dimensional.
It is evident that these representations are exactly those which we
will study from the viewpoint of dilation theory.  In this paper, we
will obtain a complete classification of these representations up to
unitary equivalence.  We will explain later how our classification
relates to the work of Bratteli and Jorgensen.

If $\fA$ is an algebra of operators, $\Lat \fA$ denotes the lattice of
all $\fA$-invar\-iant subspaces.  And if $\L$ is a lattice of
subspaces, $\Alg \L$ denotes the \wot-closed unital algebra of all
operators which leave each element of $\L$ invariant.  The algebra
$\fA$ is reflexive if it equals $\Alg \Lat \fA$.  For each reflexive
algebra, there is a quantitative measure of the distance to $\fA$
given by
\[
  \beta_\fA (T) = \sup_{L\in\L} \| P_L^\perp T P_L \| .
\]
It is easily seen that $\beta_\fA(T) \le \dist(T,\fA)$.  The algebra
is called {\it hyper-reflexive} if there is a constant $C$ such that
$\dist(T,\fA) \le C \beta_\fA(T)$.  The optimal $C$, if it is finite,
is called the distance constant for $\fA$.

The list of algebras known to be hyper-reflexive is rather short.
Arveson \cite{Arv_nest} showed that nest algebras have distance
constant 1, so that equality is achieved.  Christensen \cite{Chr}
showed that AF von Neumann algebras have distance constant at most 4.
Concerning the algebras studied in this paper, the first author
\cite{Dav_Toep} showed that the analytic Toeplitz algebra has distance
constant at most 19; and with Pitts \cite{DP1}, that all atomic free
semigroup algebras where shown to have distance constant at most 51.
The worst case for these estimates was the algebra $\fL_n$.  However a
recent general result of Bercovici \cite{Berc} applies to show that
$\fL_n$ actually has a distance constant no greater than 3.

\section{Main Results}\label{S:main}

In this paper, we generally take $n$ to be a finite integer with $n \ge2$.
However, Popescu's version of the dilation theorem is valid for $n=\infty$,
as are the results of \cite{DP1,DP2} on the structure of $\fL_n$ which
we shall use.  So the results of this paper go through for $n=\infty$
with only a few minor changes in notation, not in substance.
For ease of presentation, we will write this paper as though $n$ were
finite, and let the interested reader interpolate the $n=\infty$ case.

Consider a contractive $n$-tuple $A = (A_1,\dots,A_n)$ acting on a
finite dimensional space $\V$ of dimension $d$; i.e.
$\sum_{i=1}^n A_iA_i^* \le I$.
The Frahzo--Bunce--Popescu minimal dilation yields isometries
$S_i$ acting on a larger space $\H$.
We let $\fA$ denote the algebra generated by the $A_i$'s, and let $\fS$
be the \wot-closed algebra generated by the $S_i$'s.
We will make important use of an associated completely positive
contractive map on $\B(\V)$ given by
\[
  \Phi(X) = \sum_{i=1}^n A_i X A_i^* .
\]
The operator $\Phi^\infty(I) := \lim_{k\to\infty} \Phi^k(I)$ will also
be useful.

The first fairly easy observation is that the dilation is of Cuntz
type ($\sum_{i=1}^n S_iS_i^* = I$) if and only if
$\sum_{i=1}^n A_iA_i^* = I$; or equivalently $\Phi(I)=I$.
In general, we define the {\it pure rank} of $\fS$ to be the
multiplicity of the left regular representation in the Wold
decomposition of $\fS$.  This is the dimension of the wandering space
$\W = \ran ( I -  \sum_{i=1}^n S_iS_i^*)$.
Simple examples show that this
wandering space need not be contained in $\V$, and that even when this
pure rank is one, the pure part may have large intersection with
$\V$.  Nevertheless, it turns out that this pure rank may be easily
computed as
\[
  \prank( \fS ) = \rank( I - \Phi(I) )
              = \rank \Bigl( I - \sum_{i=1}^n A_iA_i^* \Bigr) .
\]

The irreducible summands of Cuntz type are determined by the minimal
$\fA^*$-invariant subspaces $\M$ of $\V$ on which
$\sum_{i=1}^n A_iA_i^*|_\M = I_\M$.  Such a subspace generates an
invariant subspace $\H_\M = \fS[\M]$ for $\fS$ which is necessarily
reducing. The restriction $\fS|_{\H_\M}$ of $\fS$ to this subspace is
isomorphic to the algebra $\fB_{n,m}$, where $m = \dim \M$, described
in the Background section. A crucial feature is that the projection
$P_\M$ belongs to this algebra.  This makes it possible to show that
the restriction of the $n$-tuple $A$ to $\M$ is a unitary invariant
for the dilation.

The subspace $\td{\V}$ spanned by all the minimal $\fA^*$-invariant
subspaces of this type completely determines the Cuntz part of the
dilation.  The restriction of $\fA^*$ to $\td{\V}$ is a finite
dimensional C*-algebra.  The well-known invariants for a finite
dimensional C*-algebra allow one to compute the multiplicities of each
irreducible subrepresentation.
In general, this information may be used to completely decompose the
representation into a direct sum of finitely many irreducible
representations of the types given above.  This yields complete
unitary invariants: the pure rank and the unitary equivalence class of
the restriction of $A^*$ to $\td{\V}$.

For example, one can show that $\fS$ is irreducible if and only if
either\\
(1) $\rank(I - \Phi(I)) = 1$ and $\Phi^\infty(I) = 0$, the pure
case, or\\
(2) $\{ X : \Phi(X) = X \} = \bbC I$, the Cuntz case.

The algebras $\fL_n$ and $\fB_{n,d}$ were shown to be hyper-reflexive
in \cite{DP1}.  This analysis can be used to show that all of these
algebras $\fS$ determined by a finite rank $n$-tuple are
hyper-reflexive.  The constant 51 of that paper may be improved to
5 using recent results of Bercovici \cite{Berc} which show that
the distance constant for $\fL_n$ is at most 3.

Then we turn our attention to similarity.  If two contractive
$n$-tuples are similar, it follows that their Cuntz parts are
unitarily equivalent.  However, the pure rank can change.  Indeed,
this rank can be 0 in one case and non-zero in a similar $n$-tuple.

The major interest lies in the pure case.  In this case, the algebra
$\fS$ is unitarily equivalent to a multiple of $\fL_n$, and thus
completely isometrically isomorphic and weak-$*$ homeomorphic to
$\fL_n$.  The compression $\Phi_A$ of $\fS$ to $\V$ is thus a weak-$*$
continuous representation of $\fL_n$.  The study of these
representations was initiated in \cite{DP2}.  The kernel of such a
representation is a \wot-closed ideal.  A \wot-closed ideal $\fJ$ of
$\fL_n$ is determined \cite[Theorem 2.1]{DP2} by its range
$\M = \ol{\fJ \K_n}$, which is a subspace invariant for both $\fL_n$
and its commutant $\fR_n$.  Thus we consider the associated
representation of $\fL_n$ obtained as restriction to $\M^\perp$.  This
has the same kernel $\fJ$.  In the case of an irreducible $n$-tuple,
the minimal $\fL_n^*$-invar\-iant subspaces of $\M^\perp$ yield all of
the $n$-tuples similar to $A$ which have pure rank 1.  These are the
extreme points of all such representations in the sense that $A$ can
be recovered as a C*-convex combination of them.

In particular, it follows that two irreducible $n$-tuples of matrices
are similar if and only if the induced representations of $\fL_n$ have
the same kernel.  The range space $\M$ of $\fJ= \ker \Phi_A$ has a
wandering space of dimension $1 + (n-1)d^2$.  A basis for this
wandering space yields a finite set of generators for $\fJ$ as a
\wot-closed right ideal.  They determine a corresponding finite set of
isometries $X_j$ in $\fL_n$ with the property that another contractive
$n$-tuple of $d \times d$ matrices $B$ is similar to $A$ if and only
if $\Phi_B(X_j)=0$ for $1 \le j \le 1 + (n-1)d^2$.  Since one merely
requires generators for $\fJ$ as a two-sided ideal, normally this
number of tests can be reduced. The set of isometries are canonical,
but they are generally not polynomials.  A set of polynomial
invariants can be obtained by an approximation argument.

\section{Wandering Subspaces}\label{S:wandering}

Let $\V$ be a $d$-dimensional space (possibly infinite), and let
$A_1,\dots,A_n$ be an $n$-tuple of operators in $\B(\V)$ such that
$\sum_{i=1}^n A_iA_i^* \le I$.
The Frahzo--Bunce--Popescu minimal dilation yields isometries
$S_i$ on a larger space $\H$.  Let $P_\V$ denote the
projection of $\H$ onto $\V$.   We let $\fA$ denote the algebra
generated by the $A_i$'s and $\fS$ be the \wot-closed algebra
generated by the $S_i$'s.  We first identify $\V^\perp$.

\begin{lem}\label{wander}
The subspace $\W = (\V + \sum_{i=1}^n S_i \V ) \ominus \V$ is a
wandering subspace for $S$, and
$\sum^\oplus_{v\in\F_n} S_v \W = \V^\perp$.
\end{lem}

\Prf  $\W$ is contained in $\V^\perp$, which is invariant for $S$.
Thus $S_u \W$ is orthogonal to $\V$ for every word $u\in\F_n$.
Consequently, when $|u|\ge1$, $S_u \W$ is also orthogonal to
$S_j \V$, $1 \le j \le n$.  It follows that
$S_u \W$ is orthogonal to $\V + \sum_{i=1}^n S_i \V$, which
contains $\W$.  Therefore $\W$ is wandering.
Minimality ensures that
\[
  \H = \spn \{ S_u \V : u \in \F_n \}
     = \spn \{ \V, S_u \W : u \in \F_n \} .
\]
Since $\W$ lies in the invariant subspace $\V^\perp$, this can
only occur because $\sum^\oplus_{v\in\F_n} S_v \W = \V^\perp$.
\bx

Thus $\K = \V^\perp$ is unitarily equivalent to a multiple
$\K_n^{(\alpha)}$ of Fock space, where $\alpha = \dim \W$, and
$S_i|_\K \simeq L_i^{(\alpha)}$.  Hence decomposing $\H = \V \oplus
\K$, we may write each $S_i$ as a matrix $S_i =
\begin{sbmatrix} A_i & 0 \\ X_i & L_i^{(\alpha)}\end{sbmatrix}$.

\begin{rem}\label{dim_alpha}
The range of $\sum_{i=1}^n S_iS_i^*$ includes \vspace{.2ex}
$\sum_{i=1}^n S_i \V^\perp = (\V + \W)^\perp$ as well as
$\sum_{i=1}^n S_i \V$.  Hence
$\sum_{i=1}^n S_i S_i^* = I$ if and only if
$\sum_{i=1}^n S_i \V$ contains $\V$.
Since $\V$ is invariant for $S_i^*$ and $S_i^*|_\V = A_i^*$,
\[
  \sum_{i=1}^n A_i A_i^* = \sum_{i=1}^n P_\V S_i P_\V S_i^*|_\V
  = P_\V \sum_{i=1}^n S_i S_i^*|_\V .
\]
Therefore $\sum_{i=1}^n S_i S_i^* = I$ if and only if its range
contains $\V$ if and only if $\sum_{i=1}^n A_iA_i^* = I_\V$.

Let $d = \dim \V$ be finite, and let $\alpha = \dim \W$.
Then $\alpha$ can be as large as $nd$ and as small as
$(n-1)d$.  This is easily seen since $\sum_{i=1}^n S_i \V$ is an
orthogonal direct sum and thus has dimension $nd$, so that
$\W = \bigl( \V + \sum_{i=1}^n S_i \V \bigr) \ominus \V$
can have no larger dimension than $nd$, and is at least $(n-1)d$.

When $\sum_{i=1}^n A_iA_i^* = I_\V$, we showed above that
$\sum_{i=1}^n S_iS_i^* = I$.
Then
\[
 \V = \sum_{i=1}^n S_i S_i^* \V =
      \sum_{i=1}^n S_i A_i^* \V  \subset \sum_{i=1}^n S_i \V.
\]
Hence $\W = \bigl( \sum_{i=1}^n S_i \V \bigr) \ominus \V$ has
dimension $(n-1)d$.

The case $\dim \W = nd$ occurs, for example, if $A_i = 0 $ for
$1 \le i \le n$.  The minimal dilation is just $L_i^{(d)}$.
Indeed, if $x,y \in \V$, then
\[
  (S_i x, y) = ( x, S_i^* y) =  ( x, A_i^* y) = 0
  \qfor 1 \le i \le n.
\]
Thus $\V$ is orthogonal to $\sum_{i=1}^n S_i \V$.
Therefore $\W = \sum_{i=1}^n S_i \V$ has dimension $nd$.

It is easy to combine these examples to obtain any integer in between.
\end{rem}

The $n$-tuple of isometries $S$ is called {\it pure} if it is
unitarily equivalent to a multiple of the left regular
representation.  Bunce \cite{Bun} shows that whenever $\|A\| < 1$,
the dilation $S$ is pure.  Popescu \cite{Pop_diln} shows that the
dilation is pure if and only if $\wotlim\limits_{k\to\infty}
\sum\limits_{|v|=k} A_v A_v^* = 0$.

Lemma~\ref{wander} shows that beginning with an $n$-tuple, we will
always obtain wandering vectors except when the $A_i$'s already are
isometries and $\sum_{i=1}^n A_iA_i^* = I$, in which case the dilation
is just the $A_i$'s themselves.  When there
are wandering vectors, each generates a subspace $\M$ on which the
isometries $S_i$ are unitarily equivalent to the left regular
representation.  In particular, the non-$*$ algebra $\fS$ that
they generate is very non-self-adjoint.  In fact, a
strong converse to this exists.  Recall that an algebra is {\it
reductive} if all of its invariant subspaces have invariant
(orthogonal) complements; and it is {\it transitive} if it has no
proper invariant subspaces at all.  It is an open problem equivalent
to the transitive algebra variant of the invariant subspace problem
\cite{DyerP} whether every \wot-closed reductive algebra is
self-adjoint.

\begin{lem}\label{reductive}
Let $\fS$ be a free semigroup algebra.  Then either
$\fS$ has a wandering vector or $\fS$ is reductive.
If the latter is possible, then transitive free semigroup algebras
exist.
\end{lem}

\Prf Suppose that $\fS$ has no wandering vectors.  Let $\M$ be
an invariant subspace for $\fS$.  Then $\sum_{i=1}^n S_i\M$ must equal
$\M$; for otherwise $\M \ominus \sum_{i=1}^n S_i\M$ is wandering.
Thus $\sum_{i=1}^n S_i\M^\perp = \M^\perp$, so that $\M^\perp$ is
also invariant.  Whence $\fS$ is reductive.

Now we invoke the direct integral theory for non-self-adjoint operator
algebras due to Azoff, Fong and Gilfeather \cite[Theorem~4.1]{AFG}.
Let $\fM$ be any masa in the commutant of $\fS$.
They show that $\fS$ may be decomposed with respect to $\fM$
as an integral of algebras which are transitive almost everywhere.
The isometries $S_i$ decompose as an integral of operators which are
isometries almost everywhere as well.  In particular, transitive
algebras generated by isometries with orthogonal ranges would exist.
\bx

At this point, we do not know if there are transitive free
semigroup algebras.  The motivation for suspecting that there may
be comes from the case $n=1$.  A unitary operator is reductive
unless Lebesgue measure on the whole circle is absolutely
continuous with respect to the spectral measure of the unitary
\cite{Werm}.  This is the case for `most' unitaries.

We will frequently construct {\it reducing} subspaces of $\fS$
from $\fA^*$-invar\-iant subspaces.  This procedure preserves
orthogonality as well.

\begin{lem}\label{reducing}
Suppose that $\V$ contains an $\fA^*$-invar\-iant subspace $\V_1$.
Then $\H_1 = \fS[\V_1]$ reduces $\fS$.

If $\V$ contains a pair of orthogonal $\fA^*$-invar\-iant subspaces
$\V_1$ and $\V_2$,
then $\H_j = \fS[\V_j]$ for $j=1,2$ are mutually orthogonal.

If in addition $\V = \V_1 \oplus \V_2$, then $\H$ decomposes as
$\H_1 \oplus \H_2$  and $\H_j \cap \V = \V_j$ for $j=1,2$.
\end{lem}

\Prf Since $\V_1$ is invariant for $A_i^*$, it is also invariant
for $S_i^*$.  The $\fS$-invar\-iant subspace $\H_1 = \fS[V_1]$
is spanned by vectors of the form $S_w x$ where $x \in \V_1$ and
$w\in\F_n$. Notice that $S_i^* S_w x$ equals $S_{w'} x$ if
$w = i w'$, $0$ if $w = i' w'$ for some $i' \ne i$, and $S_i^* x$
if $w = e$.  Since $\V_1$ is invariant for $\fS^*$, each of these
possibilities belongs to $\H_1$.  Thus $\H_1$ reduces $\fS$.

Likewise, if $\V_1$ and $\V_2$ are orthogonal $\fA^*$-invariant
subspaces,  it follows that $\H_1$ and $\H_2$ are orthogonal.  For if
$v_j \in \V_j$, the inner product $(S_u v_1, S_w v_2 )$ can be reduced
by cancellation of isometries until either $u$ or $w$ is the identity
element.  Then, for example when $w=e$,
\[
  (S_u v_1, v_2 ) = (v_1, S_u^* v_2) = 0
\]
by the $\fA^*$-invariance of $\V_2$ and orthogonality.

Now suppose that $\V = \V_1 \oplus \V_2$.
Since $\H_1$ contains $\V_1$ and is orthogonal to $\V_2$, it
follows that $\H_1 \cap \V = \V_1$.  Finally, $\H_1
\oplus \H_2$ is an $\fS$-reducing subspace containing $\V$, so
it is all of $\H$ by the minimality of the dilation.
\bx

\section{Finite Dimensional n-tuples}\label{S:findim}

Now let us specialize to the case when $\V$ is finite-dimensional.
In general, we can decompose the $S_i$ into a pure part and Cuntz part.
Let $\X$ be the range of $I - \sum_{i=1}^n S_iS_i^*$, which is the
wandering space for the reducing subspace
$\H_p = \sum^\oplus_{v\in\F_n} S_v \X$.
The restriction of the $S_i$ to this space yields a multiple of
the left regular representation, where the multiplicity is $\dim \X$.
We call this quantity the {\it pure rank} of the representation.  On the
complement $\H_c = \H_p^\perp$, the restrictions of $S_i$ yield a
representation of the Cuntz algebra.  Let $P_p$ and $P_c$ denote the
projections onto $\H_p$ and $\H_c$ respectively.  It is important to
note that the projection $P_p$ does not commute with $P_\V$ in
general.  So we will obtain a method of computing this pure rank directly
from the $A_i$'s.

The key technical tool in our analysis shows that $\H_c$ is
determined by $\V_c :=\H_c \cap \V$.  This is not the case for
$\H_p$. Let $R_k$ denote the projection onto
$\sum^\oplus_{|v|=k} S_v \W$, where
$\W = \bigl( \V + \sum_{i=1}^n S_i \V \bigr) \ominus \V$;
and $Q_k = \sum_{j\ge k} R_j$.  Notice that
\[
   Q_k = \sum_{|w|=k} S_w P_\V^\perp S_w^* .
\]

On any $\fS$-invar\-iant subspace $\M$ on which the
restrictions $T_i$ of $S_i$ are pure, one has for every $x\in\M$
\[
  \lim_{k\to\infty} \sum_{|w|=k} \| P_\M S_w^* x \|^2 =
  \lim_{k\to\infty} \sum_{|w|=k} \| T_w^* x \|^2 = 0  .
\]
In particular, this applies to $\H_p$ and $\V^\perp$.
While for $x \in \H_c$, one has
\[
  \sum_{|w|=k} \|S_w^* x \|^2 = \|x\|^2   \qforal k \ge 0.
\]

\begin{lem}\label{intersect1}
Suppose that $\H_1$ is a reducing subspace for $\fS$ contained
in $\H_c$.  Let $x$ be a vector such that $P_{\H_1}x \ne 0$.
Then the subspace $\M = \fS^*[x]$ contains a
vector $v$ in $\M\cap \V_c$ with $P_{\H_1}v \ne 0$.
\end{lem}

\Prf Let $P_1$ denote the projection of $\H$ onto $\H_1$.
Fix $\eps>0$; and let $x_1 = P_1 x$.
By applying the preceding remarks to both $\V^\perp$
and $\H_p$, we may choose an integer $k$ sufficiently large that
\begin{gather*}
  \sum_{|w|=k} \| P_\V^\perp S_w^* x \|^2 = \| Q_k x \|^2 < \eps^2 \\
  \sum_{|w|=k} \| P_\V^\perp S_w^* x_1 \|^2 = \| Q_k x_1 \|^2
    < \eps^2\\
\intertext{and}
  \sum_{|w|=k} \| P_p S_w^* x \|^2 = \sum_{|w|=k} \| S_w^* P_p x \|^2
   < \eps^2 .
\end{gather*}
Since $\sum_{|w|=k} S_w S_w^* P_1 = P_1$,
\begin{align*}
  \sum_{|w|=k} \| P_\V S_w^* x_1 \|^2 &=
  \sum_{|w|=k} \bigl( \| S_w^* x_1 \|^2 - \| P_\V^\perp S_w^* x_1 \|^2
      \bigr) \\
  &= \| x_1 \|^2 - \| Q_k x_1 \|^2  > \|x_1\|^2 - \eps^2 .
\end{align*}

Let $\E_1$ denote the set of words $w$ of length $k$ such that
\[
  \| P_\V S_w^* x_1 \|^2 > \eps^{-1} \|P_\V^\perp S_w^* x\|^2 .
\]
Likewise let $\E_2$ denote the set of words $w$ of length $k$ such that
\[
  \| P_\V S_w^* x_1 \|^2 > \eps^{-1} \|P_p S_w^* x\|^2 .
\]
The set $\E_1 \cap \E_2$ is relatively large in the sense that
\begin{align*}
  \sum_{w\in\E_1 \cap \E_2} \| P_\V &S_w^* x_1 \|^2 >
  \|x_1\|^2 - \eps^2 - \sum_{w\not\in\E_1} \|P_\V S_w^* x_1\|^2
          - \sum_{w\not\in\E_2} \|P_\V S_w^* x_1\|^2 \\ &>
  \|x_1\|^2 - \eps^2
      - \sum_{w\not\in\E_1} \eps^{-1} \| P_\V^\perp S_w^* x\|^2
      - \sum_{w\not\in\E_2} \eps^{-1} \| P_p S_w^* x\|^2
  \\&> \|x_1\|^2 - \eps^2 - \eps - \eps > \|x_1\|^2/4
\end{align*}
for small $\eps$.  Now we also have
$\sum_{w\in\E_1 \cap \E_2} \|P_\V S_w^* x\|^2 \le \|x\|^2 $.
Therefore there is a word $w$ in $\E_1 \cap \E_2$ such that
\[
  \|P_\V S_w^* x_1\| > \frac{\|x_1\|}{2\|x\|} \|P_\V S_w^* x\| .
\]

In this way, construct a sequence of words $w_k$ corresponding to
$\eps_k = 1/k$.  Hence define unit vectors
$y_k = S_{w_k}^* x / \|S_{w_k}^* x\|$ with the properties that
\[
  \lim_{k\to\infty} \| P_\V^\perp y_k \| \le
  \lim_{k\to\infty} \frac1{\sqrt{k}}
    \frac{ \| P_\V S_{w_k}^* P_1 x \|}{\| S_{w_k}^* x\|} =
  \lim_{k\to\infty} \frac1{\sqrt{k}}
    \frac{ \| P_\V P_1 S_{w_k}^* x \|}{\| S_{w_k}^* x\|} = 0 .
\]
Similarly,
\[
  \lim_{k\to\infty} \| P_p y_k \|  = 0 .
\]
Also
\begin{align*}
  \| P_1 y_k \| &= \| S_{w_k}^* x_1 \| / \|S_{w_k}^* x\|
  \ge \| P_\V S_{w_k}^* x_1 \| / \| S_{w_k}^* x \| \\
  &> \frac{\|x_1\|}{2\|x\|}
     \frac{ \| P_\V S_{w_k}^* x \|}{\|S_{w_k}^* x\| }
  = \frac{\|x_1\|}{2\|x\|} \| P_\V y_k \|.
\end{align*}

By the compactness of the unit ball in $\V$, there is a subsequence of
the $y_k$'s which converges to a unit vector $v$ in $\V$.
Clearly, $P_p v = 0$, and thus $v$ belongs to $\V_c \cap \fS^*[x]$;
whence this subspace is non-zero.  By construction,
$\|P_1 v\| \ge \|x_1\| / 2\|x\|$, and therefore is also non-zero.
\bx

\begin{cor}\label{intersect2}
Every non-zero subspace of $\H_c$ which is invariant for $\fS^*$
has non-zero intersection with $\V_c$.
In particular $\H_c = \fS[\V_c]$.
\end{cor}

\Prf  Let $\M$ be any non-zero $\fS^*$-invar\-iant subspace contained
in $\H_c$.  If $x$ is any non-zero vector in $\M$, the previous lemma
applied to $x$ and $\H_1 = \H_c$ shows that $\fS^*[x]$
intersects $\V_c$ non-trivially.

By Lemma~\ref{reducing}, $\N = \fS[\V_c]$ reduces $\fS$.
We claim that $\N = \H_c$.  For otherwise, let $\H_1 =
\H_c \cap \N^\perp$.  By the first paragraph, this reducing
subspace for $\fS$ must intersect $\V_c$ non-trivially.  So
$\H_1$ is not orthogonal to $\N$, contrary to fact.  Therefore
$\H_1$ must be zero.
\bx

\begin{cor}\label{intersect3}
Suppose that $\sum_{i=1}^n A_i A_i^* = I$ and $\fA = \B(\V)$.
Then every invariant subspace of $\fS^*$ contains $\V$.
\end{cor}

\Prf Since $\H = \H_c$, any $\fS^*$-invar\-iant subspace $\M$
intersects $\V$ in a non-trivial subspace.  This subspace is invariant
for $\fS^*|_\V = \fA^* = \B(\V)$.  Hence it is all of $\V$.
\bx

Let $\fB$ denote the \wot-closed operator algebra on
$\H = \V \oplus \K_n^{(\alpha)}$ spanned by
$\B(\H)P_\V$ and $0_\V \oplus \fL_n^{(\alpha)}$.

\begin{lem}\label{propA1}
Every weak-$*$ continuous functional on $\fB$ is given by a trace
class operator of rank at most $d+1$.
\end{lem}

\Prf An element $B$ of $\fB$ is determined by $BP_\V$ and
$BP_\V^\perp$.  If $e_1,\dots,e_d$ is a basis for $\V$, the former
is determined by the vectors $Be_j$.  The latter term is
unitarily equivalent to $A^{(\alpha)}$ for some $A\in \fL_n$.
Any functional $\phi$ is thus determined by a functional $\phi_0$
on $\fL_n^{(\alpha)}$ and by $d$ functionals on $\H$ given by the
Riesz Representation Theorem by a vector $y_j$.  By
\cite[Theorem~2.10]{DP1}, the functional on $\fL_n$ is given by a
rank one functional $\phi_0(A) = (A\eta,\zeta)$.  Whence
\[
  \phi(B) = \sum_{j=1}^d (Be_j,y_j) +  (B\eta,\zeta) .
\]
\upbx

\begin{cor}\label{wot=w*}
The \wot\ and weak-$*$ topologies coincide on $\fB$, and
thus also on $\fS$.  In particular, the weak-$*$ closed
algebra generated by the $S_i$'s coincides with $\fS$.
\end{cor}

\section{The Cuntz Case}\label{S:cuntz}

In this section, we specialize to the Cuntz case:
$\sum_{i=1}^n A_iA_i^* = I$ for which the isometric dilation
yields a representation of the Cuntz algebra.

\begin{eg}\label{Cuntz_state}
We begin with a description of the case in which $\V$ is one dimensional..
A special case of a finite correlated state is a Cuntz state.  This is
determined by scalars $\eta = (\eta_1,\dots,\eta_n)$
such that $\sum_{i=1}^n |\eta_i|^2 = 1$.
The state is determined by
\[
  \phi_\eta (s_{i_1}\dots s_{i_k} s^*_{j_1}\dots s^*_{j_l}) =
  \eta_{i_1}\dots \eta_{i_k} \bar{\eta}_{j_1}\dots \bar{\eta}_{j_l}.
\]
It is easy to show that the cyclic vector $\xi_\eta$ from the GNS
construction $(\H_\eta, \pi_\eta, \xi_\eta)$ spans a one-dimensional
space invariant for every $\pi_\eta(S_i^*)$.  Indeed,
\begin{align*}
  \| \pi_\eta(S_i^*) \xi_\eta - \bar{\eta}_i \xi_\eta \|^2 &=
  \lip \pi_\eta(S_i^*) \xi_\eta, \pi_\eta(S_i^*) \xi_\eta \rip -
  \eta_i \lip \pi_\eta(S_i^*) \xi_\eta, \xi_\eta \rip
   \\&\quad -
  \bar{\eta}_i \lip \xi_\eta, \pi_\eta(S_i^*) \xi_\eta \rip +
  | \eta_i |^2 \\ &=
  \phi_\eta(S_iS_i^*) - | \eta_i |^2 =
  |\eta_i|^2 - |\eta_i|^2 = 0.
\end{align*}
The restrictions $A_i^* = S_i^*|_{\spn \{ \xi_\eta \}} = \ol{\eta_i}$
satisfy $\sum_{i=1}^n A_i A_i^* = 1$.  They may be dilated to their
minimal isometric dilation, which is necessarily the original $S_i$
since $\xi_\eta$ is a cyclic vector.

Specializing to the case of $\eta = (1,0,\dots,0)$, one has $A_1 = 1$
and $A_i = 0$ for $2 \le i \le n$.  This yields the atomic representation
$\sigma_{1,1}$ mentioned in the Background section.  In particular, the
algebra $\fS$ is unitarily equivalent to $\fB_{n,1}$.

The various Cuntz states are related by the action of the gauge group
$\U(n)$ which acts as an automorphism group on $\O_n$ and on the
Cuntz--Toeplitz algebra $\E_n$.  Indeed, if we write Fock space $\K_n$
as a direct sum $\bbC \oplus \H_n \oplus \H_n^{\otimes 2} \oplus
\H_n^{\otimes 3} \oplus \dots$, where $\H_n$ is an $n$-dimensional
Hilbert space, then each unitary matrix $U\in\U(n)$ determines a
unitary operator $\td{U} = I \oplus U \oplus U^{\otimes 2} \oplus
U^{\otimes 3} \oplus \dots$ on $\K_n$.  Conjugation by $\td{U}$ acts
as an automorphism $\Theta_U$ of $\E_n$.  Moreover, it maps the ideal
of compact operators onto itself.  So it also induces an automorphism
$\theta_U$ of $\O_n$.  If $U=[u_{ij}]$ is an $n\times n$ unitary
matrix, this automorphism can also be seen to be given by
\[
 \Theta_U(L_j) = \sum_{i=1}^n u_{ij}L_i \qfor 1\le j\le n.
\]
Given $\eta$, let $U$ be any unitary with $u_{1j} = \eta_j$.  Then
it follows that
\[
  \phi_\eta(A) = \phi_{(1,0,\dots,0)}( \theta_U (A) ) \qforal A \in \O_n .
\]
So the corresponding representations are equivalent up to this
automorphism.  In particular, the {\it algebras} $\fS_\eta$ generated
by these representations are unitarily equivalent even though the
representations are not.

A crucial step in the analysis of atomic representations was to show
that certain projections lie in the algebra $\fS$.  Indeed, this is a
major advantage of $\fS$ over the C*-algebra, which contains no
non-trivial projections, and over the von Neumann algebra it
generates, which contains too many.  As a case in point, the
projection $P_\eta = \xi_\eta \xi_\eta^*$ belongs to $\fS_\eta$.
Indeed, it is the {\it only} non-trivial projection in the whole
algebra $\fS_\eta$.
\end{eg}

A crucial point of our analysis is the identification of projections
in $\fS$ in greater generality.  We begin with the irreducible case.

\begin{thm}\label{P_V}
Assume that $\sum_{i=1}^n A_iA_i^* = I$ and $\fA = \B(\V)$.
Then $\fS$ contains the projection $P_\V$.
\end{thm}

\Prf  Both $\fS$ and $P_\V$ belong to $\fB$.  If $P_\V$ were not
in $\fS$, Lemma~\ref{propA1} would provide a weak-$*$ continuous
functional $\phi$ which annihilates $\fS$ such that $\phi(P_\V)=1$.
Represent $\phi$ as a functional of rank $d+1$ in the form
$\phi(B) = \sum_{j=0}^d (Bx_j,y_j)$.
This then may be realized as a rank one functional on the
$d+1$-fold ampliation of $\fB$.  Indeed, form the vectors
$x = (x_0,\dots,x_d)$ and $y=(y_0,\dots,y_d)$.  Then
$\phi(B) = (B^{(d+1)}x,y)$.

Now the fact that $\phi$ annihilates $\fS$ means that $x$ is
orthogonal to the subspace $\M = \fS^{*(d+1)}[y]$.  The algebra
$\fS^{(d+1)}$ is generated by isometries $S_i^{(d+1)}$, which
form the minimal dilation of the $A_i^{(d+1)}$'s.  So
Corollary~\ref{intersect2} applies, and shows that $\M$
intersects $\V^{(d+1)}$ in a non-zero subspace $\M_0$ which is
invariant for $\fS^{*(d+1)}$, and thus for $\fA^{*(d+1)}$.

By hypothesis, $\fA^{*(d+1)} = \B(\V)^{(d+1)} \simeq \B(\V)
\otimes \bbC^{d+1}$, which is a finite dimensional C*-algebra.
The invariant subspace $\M_0$ is thus the range of a projection
$Q$ in the commutant $\bbC^d \otimes \fM_{d+1}$.  Let $\tilde{Q}$
denote the operator in $\bbC I_\H \otimes \fM_{d+1}$ acting on
$\H^{(d+1)}$ with the same matrix coefficients as $Q$.  That is,
$\tilde{Q}$ is the unique operator in $(\BH \otimes \bbC^{d+1})'$ such
that $P_\V^{(d+1)}\tilde{Q} = Q$.

The projection $\tilde{Q}$ yields a decomposition of $\H^{(d+1)}$ into
$\fS$-reducing subspaces $\H_1 \oplus \H_2$ where $\H_1 = \ker\tilde{Q}$
and $\H_2 = \ran\tilde{Q}$; and likewise
$\V^{(d+1)} = \V_1 \oplus \V_2$ where
\[
  \V_1 := \H_1 \cap \V^{(d+1)} = \ker Q \qand
  \V_2 := \H_2 \cap \V^{(d+1)} = \ran Q .
\]
Observe that $\M_0$ is contained in $\H_2$.  For if we had a vector
$x\in\M$ such that $P_{\H_1}x \ne 0$, then  Lemma~\ref{intersect1}
implies that there is a non-zero vector $v$ in $\M \cap \V^{(d+1)} =
\M_0$ such that $P_{\H_1}v \ne 0$.  But by definition of $Q$ and
$\td{Q}$, $\M_0$ is orthogonal to $\H_1$, a contradiction.

In particular, as $y \in \M$, we have $y = \tilde{Q}y$.  Thus
\[
  P_\V^{(d+1)} y = P_\V^{(d+1)} \tilde{Q} y = Q P_\V^{(d+1)} y
\]
belongs to $Q\V^{(d+1)} = \M_0$.
Since $x$ is orthogonal to $\M$ and hence to $\M_0$, we see that
\[
 \phi(P_\V) = ( P_\V^{(d+1)}x,y) = (x, P_\V^{(d+1)} y) = 0 .
\]
Consequently $P_\V$ belongs to $\fS$.
\bx

This immediately yields a structure theorem for $\fS$.  Note
that this does not classify the associated representations, as
they depend on the specific generators, not just the algebra.

\begin{cor}\label{P_Vcor}
Assume that $\sum_{i=1}^n A_iA_i^* = I$ and $\fA = \B(\V)$.
Then $\fS \simeq \B(\H)P_\V + \bigl( 0_\V \oplus \fL_n^{((n-1)d)} \bigr)
\simeq \fB_{n,d}$.
\end{cor}

\Prf By Theorem~\ref{P_V}, $\fS$ contains $P_\V$.  Therefore it
contains $P_\V \fS = \B(\V)$.  Moreover, it contains
$S_i P_\V^\perp \simeq 0_\V \oplus L_i^{(\alpha)}$, where
$\alpha = (n-1)d$.  Thus $\fS$ contains the \wot-closed algebra
that these operators generate, which is evidently
$0_\V \oplus \fL_n^{(\alpha)}$.
Finally, if $v$ is any non-zero vector in
$\V$, $\fS[v]$ contains $\V$ by hypothesis.  So it is all of $\H$
by minimality of the dilation.  Therefore for any $x \in \H$,
there are operators $T_k\in\fS$ such that $T_k v $ converges to
$x$.  Thus  $\fS$ contains $T_k vv^*$, which converge to the
rank one operator $xv^*$.  So $\B(\H)P_\V$ belongs to $\fS$.
This is the whole \wot-closed algebra which we called $\fB$, which
trivially contains $\fS$.  It is evident that $\fB$ is unitarily
equivalent to $\fB_{n,d}$.
\bx

Now suppose that $\fA$ is a more general subalgebra of
$\B(\V)$.  We wish to determine the structure of $\fS$ from
information about $\fA$.

\begin{lem}\label{min.cyclic}
Assume that $\sum_{i=1}^n A_iA_i^* = I$.
Suppose that $\V$ contains a minimal $\fA^*$-invar\-iant subspace $\V_0$
of dimension $d_0$ which is cyclic for $\fA$.  Then $\fS$ contains
$\B(\H)P_{\V_0}$, and is unitarily equivalent to $\fB_{n,d_0}$.
\end{lem}

\Prf  By Burnside's Theorem \cite[Corollary~8.6]{RR},
since $\fA^*|_{\V_0}$ has no proper
invariant subspaces, it must equal all of  $\B(\V_0)$.
Let $\H_0=\fS[\V_0]$.  This is a reducing subspace for
$\fS$ by Lemma~\ref{reducing}.  We will argue that $\H_0 = \H$.

Suppose that $x$ is a non-zero vector orthogonal to $\H_0$.
By Corollary~\ref{intersect2}, $\fS^*[x] \cap \V$ contains a non-zero
vector $v$.  Moreover since $\H_0^\perp$ reduces $\fS$, $v$ is
orthogonal to $\H_0$.  Therefore $\fS^*[v] = \fA^*[v]$ is an
$\fA^*$-invariant subspace orthogonal to $\V_0$.
Since $\fA \V_0 = \V$, there is an $A\in \fA$ and $v_0\in\V_0$ such
that $A v_0 = v$.  So that
\[
  \| v \|^2 = ( Av_0, v) = (v_0, A^* v) = 0 .
\]
This contradiction establishes our claim.

Now consider the compressions
$\td{A}_i = P_{\V_0}A_i|_{\V_0} = (A_i^*|_{\V_0})^*$.
Then $\sum_{i=1}^n \td{A}_i \td{A}^*_i = I_{\V_0}$
follows from the $\fA^*$-invar\-iance of $\V_0$.  Also by
hypothesis, the algebra $\td{\fA}$ generated by the $\td{A}_i$'s
is $\B(\V_0)$.  The minimal dilation of this $n$-tuple must be
precisely the restriction of $S_i$ to $\fS[\V_0]=\H$, which
is $S_i$.  So by Corollary~\ref{P_Vcor}, it follows that $\fS$ is
unitarily equivalent to $\fB_{n,d_0}$.
\bx

The following corollary is almost immediate from the structure of
$\fB_{n,d_0}$.  We point it out in order to obtain some
non-trivial consequences.

\begin{cor}\label{min}
Assume that $\sum_{i=1}^n A_iA_i^* = I$.
If $\V$ contains a subspace $\V_0$ which is cyclic for $\fA$
and is a minimal invariant subspace for $\fA^*$, then $\V_0$ is
the unique minimal $\fA^*$-invar\-iant subspace.
\end{cor}

\Prf We have $\fA^* = \fS^*|_\V$.  So by the previous
lemma, $\fA^*$ contains $P_{\V_0}\B(\V)$.  Consequently, $\V_0$
is contained in every non-zero $\fA^*$-invar\-iant subspace.
\bx

\begin{rem}
This puts constraints on which subalgebras $\fA$ of $\B(\V)$ can be
generated by $A_i$'s which satisfy $\sum_{i=1}^n A_iA_i^* = I$.
For example, the semisimple algebra of matrices of the form
$\fA_t = \begin{sbmatrix}a&0\\(b-a)t&b\end{sbmatrix}$
for $a,b$ in $\bbC$ and a fixed $t\ne0$ is similar to the
$2 \times 2$ diagonal algebra.  Note that $\fA_t$ has two
independent vectors which are cyclic for $\fA_t$ and eigenvalues
for $\fA_t^*$, namely $e_1$ and $f_2= -\bar{t}e_1 + e_2$.
By the corollary above, this cannot equal the algebra $\fA$.
Indeed, if the generators of our algebra were
$A_i = \begin{sbmatrix}a_i&0\\(b_i-a_i)t&b_i\end{sbmatrix}$,
then a computation would show that $\sum_{i=1}^n |a_i|^2 = 1$.
Likewise considering the matrix with respect to an orthonormal
basis $\{f_1,f_2\}$ would show that $\sum_{i=1}^n |b_i|^2 = 1$.
This then forces $\sum_{i=1}^n |a_i-b_i|^2 |t|^2 = 0$.
Since $t\ne0$, this forces all the $A_i$'s to be scalar, and
hence they do not generate $\fA_t$.
\end{rem}

\begin{eg}\label{1dim}
Consider a special case of the previous corollary: if $\fA$ has a
cyclic vector $e$ which is an eigenvalue for $\fA^*$. Then $\fS$
is unitarily equivalent to $\fB_{n,1}$.
The algebra $\fA$
decomposes as $\fA= \B(\V)P_e + J\fA_1 P_e^\perp$ where $P_e$ is
the orthogonal projection onto $\bbC e$, $J$ is the injection of
$\V_1=\{e\}^\perp$ into $\V$, and $\fA_1$ is a unital subalgebra
of $\B(\V_1)$.  It is easy to see that
\[
  \Lat \fA = \{ \V, JM : M \in \Lat \fA_1 \} .
\]
Hence if $\fB_1 = \Alg \Lat \fA_1$, then
\[
  \fB := \Alg \Lat \fA = \B(\V)P_e + J\fB_1 P_e^\perp .
\]
It follows that $\fA$ is reflexive if and only if $\fA_1$ is.

Thus if $\dim \V_1 >1$, there are non-reflexive examples.
For example, consider the non-reflexive algebra $\fA_1 = \{
\begin{sbmatrix}a&0\\b&a\end{sbmatrix} : a,b \in \bbC \}$.
Take $n=3$ and let
\[
 A_1 \!=\! \begin{bmatrix}
        1 &0&0\\0& 1/\sqrt{2} &0\\0&0& 1/\sqrt{2}
       \end{bmatrix} \quad
 A_2 \!=\! \begin{bmatrix}
        0&0&0\\0&0&0\\ 1/2 & 1/2 &0
       \end{bmatrix} \quad
 A_3 \!=\! \begin{bmatrix}
        0&0&0\\ 1/\sqrt{2} &0&0\\0&0&0
       \end{bmatrix}
\]
This can be seen to satisfy $\sum_{i=1}^3 A_iA_i^* = I_3$ \vspace{.2ex}
and to generate the algebra $\fA = \biggl\{
\begin{sbmatrix}c&0&0\\d&a&0\\e&b&a\end{sbmatrix} :
a,b,c,d,e \in \bbC \biggr\}$.  This is not reflexive.

Nevertheless, $\fA^*$ has a unique minimal invariant subspace, and
thus $\fS$ is unitarily equivalent to $\fB_{3,1}$, which is
hyper-reflexive.  So there is no direct correspondence between the
reflexivity of $\fA$ and $\fS$.
\end{eg}

\begin{lem}\label{irred}
Let $A=(A_1,\dots,A_n)$ be an $n$-tuple on a finite dimensional
space $\V$ such that $\sum_{i=1}^n A_iA_i^* = I$.  Let $\fA$ be
the unital algebra that they generate.  Let $S=(S_1,\dots,S_n)$ be the
minimal isometric dilation, and $\fS$ the \wot-closed algebra they
generate.  Then $\fS$ is irreducible if and only if
$\fA^*$ has a unique minimal invariant subspace $\V_0$.
\end{lem}

\Prf If $\V_0$ is unique, then it must be cyclic for $\fA$ since
$\V\ominus \fA[\V_0]$ is an invariant subspace of $\fA^*$
orthogonal to $\V_0$.  So Lemma~\ref{min.cyclic} applies.
Since $\fS$ contains $\B(\H)P_{\V_0}$, it is evidently
irreducible.

Indeed, this conclusion follows if there is {\it any} minimal
$\fA^*$-invar\-iant subspace $\V_0$ which is cyclic for $\fA$.
By Corollary~\ref{min}, $\V_0$ is
necessarily the unique minimal $\fA^*$-invar\-iant subspace.

Finally suppose that there is a minimal $\fA^*$-invar\-iant
subspace $\V_0$ which is not cyclic.  Then as in the
first paragraph, $\V\ominus \fA[\V_0]$ is an invariant subspace
of $\fA^*$ orthogonal to $\V_0$.  Let $\V_1$ be a minimal
$\fA^*$-invar\-iant subspace contained therein.  Notice that
$\fS[\V_i]$ are pairwise orthogonal reducing subspaces for
$\fS$ by Lemma~\ref{reducing}.  Hence $\H$ contains proper
reducing subspaces, and so $\fS$ is reducible.
\bx

Now we see how to deal with the case of more than one minimal
$\fA^*$-invar\-iant subspace.  In this lemma, we do not concern
ourselves with questions of uniqueness.

\begin{lem}\label{decomp}
Assume that $\sum_{i=1}^n A_iA_i^* = I$.
There is a family of minimal $\fA^*$-invar\-iant subspaces $\V_j$
of $\V$, $1 \le j \le s$, such that $\H$ decomposes into an
orthogonal direct sum of $\H_j = \fS[\V_j]$; and the
algebras $\fS|_{\H_j}$ are irreducible.
\end{lem}

\Prf This is just a matter of choosing a maximal family of
pairwise orthogonal minimal $\fA^*$-invar\-iant subspaces, say
$\V_j$ for $1 \le j \le s$.  By Lemma~\ref{reducing}, the
subspaces $\H_j = \fS[\V_j]$ are pairwise orthogonal and reducing
for $\fS$.  Moreover a direct application of the previous lemma
applied to $\H_j$ and $\V_j$ shows that $\fS|_{\H_j}$ is irreducible.
Finally we must show that $\sum_{j=1}^{\oplus s} \H_j = \H$.
Take any vector $x$ orthogonal to this sum.
By Corollary~\ref{intersect2}, $\fS^*[x]$ intersects $\V$ in a
non-zero $\fA^*$-invar\-iant subspace orthogonal to all of the
$\H_j$'s, and thus orthogonal to all of the $\V_j$'s.  This is
contrary to construction, and so yields a contradiction.
\bx

Given an $n$-tuple $A=(A_1,\dots,A_n)$ such that
$\sum_{i=1}^n A_iA_i^* = I$, let us pick a maximal family of
mutually orthogonal minimal $\fA^*$-invar\-iant subspaces  $\V_j$
of $\V$, $1 \le j \le s$; and let $P_j = P_{\V_j}$.
>From the minimality of each $\V_j$ as an $\fA^*$-invar\-iant
subspace, we know that $P_j \fA^* P_j = \B(\V_j)$.
Set $\td{\V} = \sum_{j=1}^{\oplus s} \V_j$.
Let $\td{A}_i = P_{\td{\V}} A_i |_{\td{\V}} = ( A_i^*|_{\td{\V}})^*$
be the compression of $A_i$ to $\td{\V}$; and let $\td{\fA}$ denote the
algebra they generate in $\B(\td{\V})$.

Notice that the minimal isometric dilation of $\td{A} =
(\td{A}_1,\dots,\td{A}_n)$ is precisely $S$.  It is evident that $S$
is a joint isometric dilation of $\td{A}$.  To show that it is
minimal, it suffices to show that $\fS[\td{\V}] = \H$.  But this is
established above in Lemma~\ref{decomp}.

Our goal is to show that $\td{\fA}$ is a C*-algebra.  For the
moment, let us show that it is semisimple.  Note that $\td{\fA}$
is contained in $\sum^\oplus_{1\le j \le s} \B(\V_j)$.  Moreover
the quotient map $q_j$ of compression to $\V_j$ maps $\fA$ onto
$\B(\V_j)$.  Thus the kernel of this map is a maximal ideal.
Since $\sum^\oplus q_j = \id$ is faithful, the intersection of
all maximal ideals is $\{0\}$.  Hence $\td{\fA}$ is semisimple.

Indeed, there is a minimal family $G$ so that
$\sum_{g\in G}^\oplus q_g$ is faithful.
By the Wedderburn theory, the minimal ideal
$\fA_g = \ker \sum_{h\in G\setminus \{g\}}^\oplus q_h$
is isomorphic to $\B(\V_g)$.
But this kernel will, in practice, be supported on several of
the $\V_j$'s.
This yields a partition $\td{\V} = \sum^\oplus_{g\in G} \W_g$
where $\W_g = \sum^\oplus_{j\in G_g} \V_j$ is a sum of those
$\V_j$'s equivalent to $\V_g$.  Because $\B(\V_g)$ is simple, it
follows that there is an algebra isomorphism $\sigma_j$ of
$\B(\V_g)$ onto $\B(\V_j)$ for each $j\in G_g$ such that
\[
  \td{\fA}|_{\W_g} \simeq
  \Bigl\{ \sum_{j\in G_g}\upplus \sigma_j(X) : X \in \B(\V_g) \Bigr\}.
\]
It is well-known that every isomorphism between $\B(\V_g)$ and
$\B(\V_j)$ is spatial: $\sigma_j(X) = T_j X T_j^{-1}$ for some
invertible operator $T_j$, which is unique up to a scalar multiple.

We also need to consider the unital completely positive map
$\Phi$ on $\B(\td{\V})$ given by
\[
  \Phi(X) = \sum_{i=1}^n \td{A}_i X \td{A}_i^* .
\]
Suppose that two blocks $\V_1$ and $\V_2$ are related by a
similarity as above.  Let $B_i := P_{\V_1} A_i|_{\V_1}$ and
$C_i := P_{\V_2} A_i|_{\V_2} = TB_iT^{-1}$.
Since
\[
 \sum_{i=1}^n B_iB_i^* = I_{\V_1} \qqand
 \sum_{i=1}^n C_iC_i^* = I_{\V_2},
\]
we compute that
\begin{equation*}
  I_{\V_2} = \sum_{i=1}^n (TB_iT^{-1})(TB_iT^{-1})^* =
  T\Phi_1(T^{-1} T^{*-1})T^* ,
\end{equation*}
where $\Phi_1(X) = \sum_{i=1}^n B_i X B_i^* =
 P_1\Phi(P_1XP_1)|_{\V_1}$.
Therefore
\[
  \Phi_1(T^{-1} T^{*-1}) = T^{-1} T^{*-1} .
\]

We now study this completely positive map in order to gain
information about the structure of $\td{\fA}$.

\begin{lem}\label{expn1}
Let $\Phi(X) = \sum_{i=1}^n A_i X A_i^*$ be a unital
completely positive map on $\B(\V)$, where $\V$ is finite
dimensional. If there is a non-scalar operator $X$ such that
$\Phi(X) = X$, then $\fA^* = \Alg \{ A_1^*, \dots, A_n^* \}$ has
two pairwise orthogonal minimal invariant subspaces.
\end{lem}

\Prf Since $\Phi$ is self-adjoint and unital, there is a positive
non-scalar $X$ such that $\Phi(X)=X$.  Let $\|X\|=1$ and let $\mu$
denote the smallest eigenvalue of $X$.  Then $\M = \ker (X-I)$
and $\N = \ker(X-\mu I)$ are pairwise orthogonal non-zero
subspaces. For any unit vector $x\in \M$,
\begin{align*}
  \|x\|^2  &= ( \Phi(X) x,x) = \sum_{i=1}^n ( X A_i^* x, A_i^* x ) \\
    & \le \sum_{i=1}^n ( A_i^* x, A_i^* x ) = \|x\|^2
\end{align*}
This equality can only hold if each $A_i^* x$ belongs to $\M$.
Hence $\M$ is invariant for $\fA^*$.

This argument worked because $1$ is an extreme point in the spectrum
of $X$.  This is also the case for $\mu$.  Hence a similar argument
shows that $\N$ is invariant for $\fA^*$.
\bx

The following is a partial converse to the previous lemma.

\begin{lem}\label{expn2}
Let $\Phi(X) = \sum_{i=1}^n A_i X A_i^*$ be a unital completely
positive map on $\B(\V)$, where $\V$ is finite dimensional.
Suppose that $A_i = B_i \oplus C_i$ with respect to an orthogonal
decomposition $\V=\V_1\oplus \V_2$.  Moreover, suppose that
$\Alg\{B_i\} = \B(\V_1)$ and $\Alg\{C_i\} = \B(\V_2)$.
If there is an operator $X$ such that $\Phi(X)=X$ and
$X_{21} :=P_{\V_2} X P_{\V_1} \ne 0$, then there is a unitary
operator $W$ such that $C_i = W^*B_iW$.  Moreover the fixed
point set of $\Phi$ consists of all matrices of the form
$\begin{sbmatrix}
 a_{11}I_{\V_1}&a_{12}W^*\\a_{21}W&a_{22}I_{\V_2}
\end{sbmatrix}$.
\end{lem}

\Prf Since $\Phi$ is self-adjoint, we may suppose that $X=X^*$.
Then normalize so that $\|X_{21}\|=1$.
Let $\M = \{  v \in \V_1 : \| X_{21} v \| = \|v\| \}$.
Also let $\N = X_{21} \M$ denote the corresponding subspace of $\V_2$.
Write $B = \begin{bmatrix} B_1& \dots&B_n \end{bmatrix}$ and
$C = \begin{bmatrix} C_1& \dots& C_n \end{bmatrix}$, so that
\[
 X_{21}v = \Phi(X_{21})v = C X_{21}^{(n)} B^*v \qfor v \in \M .
\]
Since $C$ and $B^*$ are contractions, and $X_{21}$ achieves its
norm on $v$, it follows that $B^*v$ belongs to the subspace
$\M^{(n)}$ on which $X_{21}^{(n)}$ achieves its norm.
Consequently each $B_i^*$ leaves $\M$ invariant.
But as $\Alg\{B_i\} = \B(\V_1)$, this forces $\M=\V_1$.
Similarly, consideration of $X_{12}=X_{21}^*$ shows that $\N =
\V_2$. Thus $X_{21}$ and $X_{21}^*$ are isometries; so
$W = X_{21}|_{\V_1}$ is a unitary map from $\V_1$ onto $\V_2$.

Further, the identity above now shows that $W = C W^{(n)} B^*$.
Hence for all $v \in \V_1$
\[
  \| v \| = \| Wv \| = \| C W^{(n)} B^* v \|
          \le \| W^{(n)} B^* v \| \le \|v\| .
\]
In particular, $C$ acts as an isometry from the range of $W^{(n)} B^*$
onto the range $\ran W = \V_2$.  Since $C$ is contractive, it must be
zero on the orthogonal complement of $\ran W^{(n)} B^*$.  This implies
that $C^*$ is an isometry of $\V_2$ onto $\ran W^{(n)} B^*$.
Consequently, $C^* W = W^{(n)} B^*$; or equivalently,
$C_i^* = WB_i^*W^*$ for $1 \le i \le n$.


Finally, if $Y \in \B(\V_1,\V_2)$ and
$\begin{sbmatrix}0&0\\Y&0\end{sbmatrix}$ is fixed by $\Phi$, then
\[
  Y = \sum_{i=1}^n C_i Y B_i^* = \sum_{i=1}^n W B_i W^* Y B_i^*
    = W \Phi_1 (W^*Y)
\]
where $\Phi_1(X) = \sum_{i=1}^n B_i X B_i^*$ acts on $\B(\V_1)$.
By Lemma~\ref{expn1}, $W^*Y$ is scalar; so $Y$ is a multiple
of $W$.  A similar analysis works for the other coordinates.
\bx

\begin{eg}
Let
\[
 A_1 = \begin{bmatrix} 1/\sqrt{2} &0&0\\ 1/2\sqrt{2} &
 1/2 & 1/2\sqrt{2} \\ 0&0& 1/\sqrt{2} \end{bmatrix}
 \;\;\text{and}\;\;
 A_2 = \begin{bmatrix} 1/\sqrt{2} &0&0\\ -1/2\sqrt{2} &
 1/2 & -1/2\sqrt{2} \\ 0&0& 1/\sqrt{2} \end{bmatrix} .
\]
Then the matrix
$X = \begin{sbmatrix}1&0&0\\0& 1/2 &0\\0&0&0 \end{sbmatrix}$
satisfies $\Phi(X)=X$.
A calculation shows that the fixed point set of $\Phi$ is the set
of matrices $X=[x_{ij}]$ such that $x_{12}=x_{21}=x_{23}=x_{32}=0$
and $x_{11}+x_{13}+x_{31}+x_{33}=2x_{22}$.
In particular, this is not an algebra.
The algebra $\fA^*$ has two minimal invariant subspaces, $\bbC e_1$
and $\bbC e_3$.  Note that the compression of $\fA$ to
$\spn \{ e_1, e_3 \}$ consists of scalar matrices, and the fixed
point set of the restricted completely positive map is the full
$2 \times 2$ matrix algebra.
\end{eg}

We can now utilize the detailed information about the map $\Phi$
to determine the algebra $\td{\fA}$.

\begin{thm}\label{Calg}
Let $\Phi(X) = \sum_{i=1}^n A_i X A_i^*$ be a unital completely
positive map on $\B(\V)$, where $\V$ is finite dimensional.
Suppose that $\V$ is the orthogonal direct sum of minimal
$\fA^*$-invar\-iant subspaces.  Then $\fA$ is a C*-algebra and the
fixed point set of $\Phi$ coincides with the commutant of $\fA$.
\end{thm}

\Prf Let $\V = \sum^\oplus_j \V_j$ be an orthogonal decomposition
into minimal $\fA^*$-invar\-iant subspaces.  The restriction of
$\fA$ to $\V_j$ is all of $\B(\V_j)$ by Burnside's Theorem.
Thus the restriction of $\Phi$ to $\B(\V_j)$ maps onto the
scalars by Lemma~\ref{expn1}.  By the earlier analysis, $\fA$ splits
into an algebraic direct sum of minimal ideals which are isomorphic to
full matrix algebras.  These are determined by certain spatial
intertwining relations between some of the summands. If the
restriction of $A_i$'s to $\V_1$ and $\V_2$ are related by an
intertwining operator $T$, then we showed that
$\Phi_1(T^{-1}T^{*-1}) = T^{-1}T^{*-1}$.  But this is scalar by
Lemma~\ref{expn1}.  So after scaling $T$, it becomes a unitary.
It follows that $\fA$ is a C*-algebra.

Evidently, $\Phi$ fixes the commutant of $\fA=\fA^*$.  Suppose
that $\Phi(X)=X$.  If $\V_k$ and $\V_l$ are not related by a
unitary intertwining map, then by Lemma~\ref{expn2},
$P_k X P_l = 0$.  While if they are related by a unitary
$W_{kl}$, then $P_k X P_l = x_{kl} W_{kl}$ belongs to $\fA'$.  It
follows that the  fixed point set is precisely the commutant of $\fA$.
\bx

Now it is possible to provide a complete description of the
algebra $\fS$ in the Cuntz case.

\begin{lem}\label{projections}
Let $P_g$ for $g \in G$ denote the minimal central projections
in $\td{\fA}$.  These projections belong to $\fS$.
\end{lem}

\Prf We follow the lines of Theorem~\ref{P_V}.
We may work in the algebra
$\fB = \B(\H)P_\V + \bigl( 0_\V \oplus \fL_n^{(\alpha)} \bigr)$
which contains $\fS$ and each projection $P_g$.
If a central projection $P$ of $\td{\fA}$ were not in $\fS$,
by Lemma~\ref{propA1} it could be separated from $\fS$ by a
functional of rank $d+1$, which as before we write as
$\phi(A) = (A^{(d+1)} x, y)$.  Let
$\M = \fS^{*(d+1)}[y]$ and $\M_0 = \td{\V}^{(d+1)} \cap \M$.
This subspace $\M_0$ is invariant for the C*-algebra
$\td{\fA}^{*(d+1)} = \td{\fA}^{(d+1)}$, and thus is the range of a
projection $Q$ in its commutant.

Now $P^{(d+1)}$ lies in the centre of $\td{\fA}^{(d+1)}$, and
thus commutes with $Q$ as well.  Therefore $\td{\V}^{(d+1)}$
decomposes as
\begin{gather*}
  P^{(d+1)} Q \td{\V}^{(d+1)} \oplus
  P^{\perp(d+1)} Q \td{\V}^{(d+1)} \oplus
  P^{(d+1)} Q^\perp \td{\V}^{(d+1)} \oplus
  P^{\perp(d+1)} Q^\perp \td{\V}^{(d+1)} \\ =:
  \M_{pq} \oplus \M_{p^\perp q} \oplus \M_{pq^\perp} \oplus
  \M_{p^\perp q^\perp} .
\end{gather*}
This determines an orthogonal decomposition of  $\td{\V}^{(d+1)}$
into four reducing subspaces for $\td{\fA}^{(d+1)}$.  Note that $\M_0$
is the sum of the first two.  Recall the remarks following
Lemma~\ref{decomp} that $S$ is the minimal isometric dilation of
$\td{A}$.  So by Lemma~\ref{reducing},
$\H^{(d+1)}$ has an orthogonal decomposition into the four
reducing subspaces for $\fS^{(d+1)}$ generated by these
subspaces of $\td{\V}^{(d+1)}$, say
\[
  \H^{(d+1)} = \H_{pq} \oplus \H_{p^\perp q} \oplus
               \H_{pq^\perp} \oplus \H_{p^\perp q^\perp} .
\]
Moreover, Lemma~\ref{intersect1} shows as in the proof of
Theorem~\ref{P_V} that $y$ belongs to
$\H_{pq} \oplus \H_{p^\perp q} = \fS^{(d+1)}[\M_0]$.

It is evident from this construction that each of these four
subspaces $\H_{ij}$ is mapped onto the corresponding $\M_{ij}$ by
the orthogonal projection $P_{\td{\V}}^{(d+1)}$ onto
$\td{\V}^{(d+1)}$.  Therefore, since $P^{(d+1)}$ is dominated by
this projection, it is clear that it maps $y$ into $\M_{pq}$,
which is contained in $\M_0$. As before, we obtain that $x$ is
orthogonal to $\M_0$, and therefore $\phi(P) = 0$.
Hence we conclude that $P$ belongs to $\fS$.
\bx

\begin{thm}\label{Cuntz_case}
Let $A_1,\dots,A_n$ be operators on a finite dimensional space
$\V$ such that $\sum\limits_{i=1}^n A_iA_i^* = I$,\vspace{.3ex}
and let $S_1,\dots,S_n$ be their joint isometric dilation.
Let $\td{\V}$ be the subspace of $\V$ spanned by all minimal
$\fA^*$-invar\-iant subspaces.  Then the compression $\td{\fA}$ of
$\fA$ to $\td{\V}$ is a C*-algebra.
Let $\td{\fA}$ be decomposed as
$\sum^\oplus_{g\in G} \fM_{d_g} \otimes \bbC^{m_g}$
with respect to a decomposition
$\td{\V} = \sum^\oplus_{g\in G} \V_g^{(m_g)}$,
where $\V_g$ has dimension $d_g$ and multiplicity $m_g$.
Let $P_g$ denote the projection onto $\V_g$.
Then the dilation acts on the space
\[
  \H = \sum_{g\in G}\upplus \H_g^{(m_g)}
     = \td{\V} \oplus \K_n^{(\alpha)}
\]
where $\H_g = \V_g \oplus \K_n^{(\alpha_g)}$ and
$\alpha_g = d_g(n-1)$ and
\[
\alpha = \sum_{g\in G} \alpha_g m_g = (n-1) \sum_{g\in G} d_g m_g .
\]
The algebra $\fS$ decomposes as
\[
  \fS \simeq
      \sum_{g\in G}\upplus \bigl( \B(\H_g)P_g \bigr)^{(m_g)}
      \;+\; \bigl( 0_{\td{\V}} \oplus \fL_n^{(\alpha)} \bigr) .
\]
\end{thm}

\Prf This is now just a matter of putting the pieces together and
clearing up some final details.  Let $\V_g$, $1 \le g \le s$, be
any maximal family of pairwise orthogonal minimal
$\fA^*$-invar\-iant subspaces.
Let $\td{\V} = \sum^\oplus_{1 \le g \le s} \V_g$.
(Do not worry at this stage about the uniqueness of the
definition of $\td{\V}$.)
By Lemma~\ref{reducing}, $\H_g = \fS[\V_g]$ are pairwise
orthogonal reducing subspaces of $\fS$.
Let $\M =  \sum^\oplus_{1 \le g \le s} \H_g$.
We claim that $\M = \H$.  Indeed, were there a non-zero vector in
$\M^\perp$, then by  Corollary~\ref{intersect2}, $\M^\perp \cap \V$
would be an non-zero $\fA^*$-invar\-iant subspace orthogonal to
$\td{\V}$, contrary to fact.

It now follows as above that if we compress each $A_i$ to $\td{A}_i$ on
$\td{\V}$, then this new $n$-tuple has the identical joint
isometric dilation $S_i$, and it is the minimal dilation by the
previous paragraph.
By Theorem~\ref{Calg}, the algebra $\td{\fA}$ that they generate
is self-adjoint.  Then applying Lemma~\ref{projections}, we
deduce that the projection onto $\td{\V}$ belongs to $\fS$, and
that $P_{\td{\V}} \fS = \td{\fA}$.

The restriction of $\fS$ to each reducing subspace $\H_g$ is
isomorphic to $\fB_{n,d_g}$.  Moreover the
restriction of $\fS$ to $\td{\V}^\perp$ is canonically isomorphic
to $\fL_n^{(\alpha)}$, where by canonical we mean that
$u(S)|_{\td{\V}^\perp} \simeq L_u^{(\alpha)}$ when we make the
natural identification of ${\td{\V}^\perp}$ with $\K_n^{(\alpha)}$
as in Lemma~\ref{wander}.

Now the finite dimensional C*-algebra $\td{\fA}$ may be
decomposed as $\sum^\oplus_{g\in G} \fM_{d_g} \otimes \bbC^{m_g}$.
The multiplicities reflect the fact that the restrictions of
$A_i^*$ to different $\V_g$'s may be unitarily equivalent.  As
before, choose a maximal subset $G$ of pairwise inequivalent
subspaces $\V_g$, and let $\W_g = \sum^\oplus_{j\in G_g} \V_j$ be
the sum of all subspaces equivalent to $\V_g$.  Then $\W_g$ may be
naturally identified with $\V_g \otimes \bbC^{m_g}$ so that
$A_i^*|_{\W_g} \simeq \bigl( A_i^*|_{\V_g} \bigr)^{(m_g)}$.
This identifies $\td{\V}$ with $\sum^\oplus_{g\in G} \V_g^{(m_g)}$.

By the uniqueness of the minimal isometric dilation, it also
follows that there is a corresponding unitary equivalence between
$\sum^\oplus_{j\in G_g} \H_j$ and $\H_g \otimes \bbC^{(m_g)}$ so
that the restriction of $S_i$ is identified with
$\bigl( S_i|_{\H_g} \bigr)^{(m_g)}$.  By Lemma~\ref{projections},
the projection $P_{\W_g} \simeq P_g^{(m_g)}$ belongs to $\fS$.
Thus we now see that $\fS P_{\td{\V}}$ decomposes as
$\sum^\oplus_{g\in G} \bigl( \B(\H_g)P_g \bigr)^{(m_g)}$.
Combining all of the pieces, we obtain the desired structure
theory for $\fS$.

It remains to establish the uniqueness of $\td{\V}$.
We can now see that $P_{\td{\V}}$ is the unique maximal finite
rank projection in $\fS$.
Indeed, every operator in $\fS$ has a lower triangular form with
respect to the decomposition $\H = \td{\V} \oplus \K_n^{(\alpha)}$.
By \cite[Corollary~1.8]{DP1}, $\fL_n$ contains no proper
projections.  Therefore all finite rank projections are supported
by $\td{\V}$.
Now suppose that $\V_0$ is any minimal $\fA^*$-invar\-iant
subspace.  It may be extended to a maximal family of pairwise
orthogonal minimal $\fA^*$-invar\-iant subspaces, and the
construction may proceed as above.  The same subspace $\td{\V}$
necessarily is obtained by the uniqueness of this maximal
projection. In particular, $\td{\V}$ must contain {\it every}
minimal $\fA^*$-invar\-iant subspace.  Thus it is the span of all
such subspaces.
\bx

\section{The General Finite Dimensional Case}\label{S:general}

We now return to the problem posed in Section~\ref{S:findim}.
Starting with a contractive $n$-tuple $A_1,\dots,A_n$ with
minimal joint isometric dilation $S_1,\dots,S_n$, we wish to
understand the structure of $\fS = \Alg\{S_1,\dots,S_n\}$ in
terms of the structure of the $n$-tuple $A$ and the algebra $\fA$
that it generates.

Recall from the discussion in Section~\ref{S:findim} that
$\H = \H_p \oplus \H_c$, where $\H_p$ is the pure part determined
by the wandering subspace of $S$, and $\H_c$ is the Cuntz part;
and that $P_p$ and $P_c$ denote the orthogonal projections
onto these subspaces.  We need a method of getting information
about this decomposition from $A$.  Corollary~\ref{intersect2}
shows that $\H_c = \fS[\V_c]$, so $\H_c$ is recovered if we can
compute $\V_c$.

Again we consider the completely positive map $\Phi(X) =
\sum_{i=1}^n A_i X A_i^*$.  This is no longer unital, since
$\Phi(I) = AA^* = \sum_{i=1}^n A_i A_i^* \le I$.  But it is
completely contractive.  Thus the sequence $\Phi^k(I)$ is a
decreasing sequence of positive operators, and therefore converges to
a limit which we denote as $\Phi^\infty(I)$.

\begin{lem}\label{expn_proj}
$\Phi^\infty(I) = P_\V P_c P_\V$.
Hence $\V_c = \ker(I-\Phi^\infty(I))$.
\end{lem}

\Prf If $x\in\H_c$,
\[
  \sum_{|w|=k} \| S_w^* x \|^2 = \|x\|^2  .
\]
On the other hand, any vector $x$ in $\H_p$ satisfies
\[
  \lim_{k\to\infty} \sum_{|w|=k} \| S_w^* x \|^2 = 0 .
\]
Thus if $x$ is any vector in $\H = \H_c \oplus \H_p$,
\[
  \lim_{k\to\infty} \sum_{|w|=k} \| S_w^* x \|^2 = \| P_c x \|^2 .
\]
We write $A_w^* := w(A)^* = S_w^*|_{\V}$.  Now if $v \in \V$,
\[
  \lim_{k\to\infty} \sum_{|w|=k} \| A_w^* v \|^2 =
  \lim_{k\to\infty} \sum_{|w|=k} \| S_w^* v \|^2 = \| P_c v \|^2 .
\]

It is evident that $\Phi^k(I) = \sum_{|w|=k}A_wA_w^*$
and thus
\[
  (\Phi^k(I) v,v) = \sum_{|w|=k} \| A_w^* v \|^2 .
\]
Therefore
\[
  (\Phi^\infty(I) v,v) = \| P_c v\|^2 = (P_\V P_c P_\V v,v).
\]
Since a sesquilinear form can be recovered from its quadratic
form by the polarization identity, it follows that
$\Phi^\infty(I) = P_\V P_c P_\V$.

In particular, $\ker(I-\Phi^\infty(I)) = \V \cap \H_c = \V_c$.
\bx

We have $\H_c = \fS[\V_c]$, and thus the restriction of the
$S_i$'s to $\H_c$ are the minimal joint isometric dilations of
the compressions of the $A_i$'s to $\V_c$.
By the previous section, we know that $\fS|_{\H_c}$ is determined
by the restriction of $\fA$ to the span $\td{\V}$ of all
$\fA^*$-invar\-iant subspaces contained in $\V_c$.  It is
desirable to give a definition that is somewhat independent of the
definition of $\V_c$.
{\it the space $\td{\V}$ is the span of all $\fA^*$-invar\-iant
subspaces $\W$ on which $\sum_{i=1}^n A_iA_i^*|_\W = I_\W$.}
Indeed, the condition
\[
 I_\W = \sum_{i=1}^n A_iA_i^*|_\W = \sum_{i=1}^n S_iS_i^*|_\W
\]
implies that $\W$ is contained in $\H_c$, whence in $\H_c \cap \V
= \V_c$.  Thus $\W$ is contained in $\td{\V}$ by
Theorem~\ref{Cuntz_case}.  The converse follows from the
description there of $\td{\V}$.

\begin{lem}\label{Cuntz_projection}
The projection $P_{\td{\V}}$ belongs to $\fS$.
\end{lem}

\Prf We may assume that $P_{\td{\V}}\ne0$.
Suppose to the contrary that $\phi$ is a \wot-continuous
functional which separates $P_{\td{\V}}$ from $\fS$.  Then as
before, we represent $\phi(X) = (X^{(d+1)}x,y)$ on an algebra
$\fB$ containing $\fS$ and $P_{\td{\V}}$.

Split $x = x_c \oplus x_p$ and $y = y_c \oplus y_p$ corresponding
to the decomposition of
$\H^{(d+1)} = \H_c^{(d+1)} \oplus H_p^{(d+1)}$.
The functional $\phi_p(X) = (X^{(d+1)} x_p, y_p)$ acts on the
pure part $\fS|_{\H_p}$.  Since the Cuntz part is non-zero and
contains wandering subspaces on which the $S_i$'s are unitarily
equivalent to $L_i$, it is easy to find vectors $x_0$ and
$y_0$ in $\H_c$ such that $\phi_p(X) = (X x_0, y_0)$.  Thus
if we set $x' = x_c \oplus x_0$ and $y' = y_c \oplus y_0$ in
$\H^{(d+2)}$, we obtain vectors in $\H_c^{(d+2)}$ such that
$\phi(X) = ( X^{(d+2)} x', y')$.

Thus by restricting to $\H_c$, we obtain a \wot-continuous linear
functional which separates $P_{\td{\V}}$ from $\fS|_{\H_c}$.
This contradicts Lemma~\ref{projections}.
Thus $P_{\td{\V}}$ must belong to $\fS$.
\bx

Next we wish to compute the pure rank of the dilation.

\begin{lem}\label{pure rank}
The pure rank of $\fS$ is computed as
\[
  \prank( \fS ) = \rank( I - \Phi(I) )
              = \rank \Bigl( I - \sum_{i=1}^n A_iA_i^* \Bigr) .
\]
\end{lem}

\Prf  The wandering space is
$\X = \ran \bigl( I -  \sum_{i=1}^n S_iS_i^* \bigr)$
and the pure rank of $\fS$ equals $\dim \X$.
The minimality of the dilation means that $\X$ does not intersect
$\V^\perp$.
Therefore $P_\V P_\X P_\V$ has the same rank as $P_\X$.
However it is easy to see that
\begin{align*}
  P_\V P_\X P_\V |_\V
  &= P_\V \Bigl( I_\H - \sum_{i=1}^n S_iS_i^* \Bigr) P_\V |_\V \\
  &=  I_\V - \sum_{i=1}^n A_iA_i^* = I_\V - \Phi(I_\V) .
\end{align*}
Thus $\prank(\fS) = \rank \bigl( I - \Phi(I) \bigr)$.
\bx

\begin{eg}
Any subtlety of the preceding lemma is due to fact that $\X$ is
not, in general, contained in $\V$.  To illustrate this, consider
the following example.  Let
\[
A_1 = \begin{bmatrix}1&0&0\\0&0&0\\0&0&0\end{bmatrix} \qqand
A_2 = \begin{bmatrix}0&0&0\\1/2&0&1/2\\0&0&0\end{bmatrix} .
\]
Then
\[
 A_1A_1^* + A_2A_2^* =
\begin{bmatrix}1&0&0\\0&1/2&0\\0&0&0\end{bmatrix} .
\]
It is clear that $\bbC e_1$ and $\bbC e_3$ are pairwise orthogonal
minimal $\fA^*$-invar\-iant subspaces.  The vector $e_1$ generates the
subspace
$\H_1 = \ol{\fS e_1}$ on which the representation is equivalent to
the atomic representation $\sigma_{1,1}$.  Furthermore, $e_3$ is a
wandering vector generating a copy of the left regular representation
on $\H_3 = \ol{\fS e_3 }$.  However $e_2$ is not orthogonal to
$\H_1 \oplus \H_3$.  One can show that there is a second wandering
vector $\zeta := e_2 - P_\V^\perp S_2 (e_1 + e_3)$.  The subspace
$\H_2 = \ol{ \fS \zeta }$ yields the decomposition
$\H = \H_1 \oplus \H_2 \oplus \H_3$.

The point here is that this decomposition does not decompose $\V$
into orthogonal pieces.  In fact, $\H_2$ has trivial intersection
with $\V$; and the vector $e_2$ has components in all three pieces.
\end{eg}

We can now completely describe the algebra $\fS$ determined
by the joint isometric dilation of a contractive $n$-tuple.
There is nothing to do except combine the information in
Theorem~\ref{Cuntz_case} with the preceding two lemmas.

\begin{thm}\label{general}
Let $A_1,\dots,A_n$ be a contractive $n$-tuple on a finite
dimensional space $\V$ with joint minimal isometric dilation
$S_1,\dots,S_n$ on $\H$.  The space $\H$ decomposes as $\H_p
\oplus \H_c$ into its pure and Cuntz parts.  The multiplicity of $\H_p$
is $\prank(\fS) = \rank \bigl( I - \sum_{i=1}^n A_iA_i^* \bigr)$.
The subspace $\td{\V}$ spanned by all minimal $\fA^*$-invar\-iant
subspaces $\W$ on which $\sum_{i=1}^n A_iA_i^*|_\W = I_\W$
determines $\H_c = \fS[\td{\V}]$.

The compression $\td{\fA}$ of $\fA$ to $\td{\V}$ is a C*-algebra.
Let $\td{\fA}$ be decomposed as
$\sum^\oplus_{g\in G} \fM_{d_g} \otimes \bbC^{m_g}$
with respect to a decomposition
$\td{\V} = \sum^\oplus_{g\in G} \V_g^{(m_g)}$,
where $\V_g$ has dimension $d_g$ and multiplicity $m_g$;
and let $P_g$ denote the projection onto $\V_g$.
Then the dilation acts on the space
\[
  \H = \sum_{g\in G}\upplus \H_g^{(m_g)} \oplus \H_p
       = \td{\V} \oplus \K_n^{(\alpha)}
\]
where $\H_g = \V_g \oplus \K_n^{(\alpha_g)}$,
$\alpha_g = d_g(n-1)$ and
\begin{align*}
  \alpha &= \sum_{g\in G} \alpha_g m_g + \prank(\fS)
       \\&= (n-1) \sum_{g\in G} d_g m_g  +
               \rank \Bigl( I - \sum_{i=1}^n A_iA_i^* \Bigr).
\end{align*}
The algebra $\fS$ decomposes as
\[
  \fS \simeq
      \sum_{g\in G}\upplus \bigl( \B(\H_g)P_g \bigr)^{(m_g)}
      \;+\; \bigl( 0_{\td{\V}} \oplus \fL_n^{(\alpha)} \bigr) .
\]
\end{thm}

We now collect some of the consequences of this theorem.
First we obtain simple conditions to determine when the dilation of
$A$ is irreducible.

\begin{cor}\label{irreducible}
The algebra $\fS$ determined by the joint isometric dilation of a
contractive $n$-tuple $A$ on a finite dimensional space $\V$ is
irreducible if and only if either
\begin{itemize}
 \item[$(1)$]
   $\ran (I - \sum_{i=1}^n A_iA_i^* ) = \bbC v \ne 0$ and
   $v$ is cyclic for $\fA$.
   In this case, $\fS$ is unitarily equivalent to $\fL_n$.
 \item[{\it or}]
 \item[$(2)$]
   $\sum_{i=1}^n A_iA_i^* = I$ and $\fA^*$ has a minimal
   invariant subspace $\V_0$ which is cyclic for $\fA$.  In this
   case, $\fS$ is unitarily equivalent to $\fB_{n,d_0}$
   where $d_0 = \dim \V_0$.
\end{itemize}
which are respectively equivalent to
\begin{itemize}
 \item[$(1')$]
   $\rank(I - \Phi(I)) = 1$ and $\Phi^\infty(I) = 0$.
 \item[{\it or}]
 \item[$(2')$]
   $\{ X : \Phi(X) = X \} = \bbC I$.
\end{itemize}
\end{cor}

\Prf $\fS$ is irreducible if and only if either it is pure with
pure rank 1, or it has pure rank 0 and, by Lemma~\ref{irred}, has a
unique minimal $\fA^*$-invar\-iant subspace.

By Lemma~\ref{pure rank}, the pure rank is 1 precisely when
$\rank(I - \Phi(I)) = 1$, or equivalently that
$\ran (I - \sum_{i=1}^n A_iA_i^* )$ is a one-dimensional subspace
$\bbC v$.  Now $\fS$ is pure precisely when $\H_c=\{0\}$, which
by Corollary~\ref{intersect2} is equivalent to $\V_c = \{0\}$.
By Lemma~\ref{expn_proj}, this is equivalent to
$\Phi^\infty(I) = 0$, which establishes the equivalence with
$(1')$.  Now $\V_c$ is $\fA^*$-invar\-iant and orthogonal to $v$,
and therefore orthogonal to $\fA v$.  So if $v$ is $\fA$--cyclic,
then $\fA[v] = \V$ and $\V_c = \{0\}$.  Conversely, if $\fA[v]$ is
proper, then $\M = \fA[v]^\perp$ is $\fA^*$-invar\-iant.  But
$\sum_{i=1}^n A_iA_i^*|_\M = I_\M$ because of the condition on
$\Phi(I)$.  So $\M$ is contained in $\H_c$.  This verifies the
equivalence with (1).

The Cuntz case is synonymous with the condition
$\sum_{i=1}^n A_iA_i^* = I$.  If $\M$ is a minimal
$\fA^*$-invar\-iant subspace, then $\fA[\M]^\perp$ contains
another.  So if $\M$ is unique, it must be cyclic.  Conversely,
if it is not unique, then by Theorem~\ref{Cuntz_case}, $\td{\V}$
contains at least two pairwise orthogonal minimal
$\fA^*$-invar\-iant subspaces, one of which may be taken to be
$\M$; call the other $\M'$.  Then $\fA[\M]$ is orthogonal to
$\M'$ and thus it is not cyclic for $\fA$.  This establishes the
equivalence with (2).

Condition $(2')$ contains the fact that $\Phi(I)=I$, so this is the
Cuntz case.  If there were more than one minimal $\fA^*$-invar\-iant
subspace, then by Theorem~\ref{Calg} the fixed point
algebra contains non-scalar operators.  Conversely, if $\Phi$ has
non-scalar fixed points, then Lemma~\ref{expn1} shows that there
are two orthogonal $\fA^*$-invar\-iant subspaces.  So $(2')$ is
equivalent to irreducibility.
\bx

\begin{cor}\label{pure}
The minimal isometric dilation of a finite dimensional $n$-tuple
$A = (A_1,\dots,A_n)$ is pure if and only if
$\fA (I - \sum_{i=1}^n A_iA_i^*) \V = \V$ or equivalently that
$\Phi^\infty(I) = 0$.
\end{cor}

\Prf The dilation has a Cuntz part if and only if there is a
$\fA^*$-invariant subspace $\M$ contained in
$\ker (I - \sum_{i=1}^n A_iA_i^*)$.  This is equivalent to having the
proper $\fA$-invariant subspace $\M^\perp$ containing
\[
  \bigl( \ker (I - \sum_{i=1}^n A_iA_i^*) \bigr)^\perp =
  \ran (I - \sum_{i=1}^n A_iA_i^*) .
\]
The minimal such subspace is
clearly $\fA (I - \sum_{i=1}^n A_iA_i^*) \V$.  Thus the dilation is
pure precisely when $\fA (I - \sum_{i=1}^n A_iA_i^*) \V = \V$.

Evidently, if there is a Cuntz part, then
\[
  \Phi^\infty(I) \ge \Phi^\infty( P_{\td{\V}} ) = P_{\td{\V}} .
\]
Conversely, if $A$ is pure, then
$\sotlim_{k\to\infty} \sum_{|w|=k} S_wS_w^* = 0$.  The compression of
$S_wS_w^*$ to $\V$ is $A_wA_w^*$, and thus
\[
  \sum_{|w|=k} P_\V S_wS_w^*|_\V = \sum_{|w|=k} A_wA_w^* = \Phi^k(I) .
\]
Since $\V$ is finite dimensional, this converges to 0 in norm.
\bx

Our theorem also provides simple complete unitary invariants for the
associated finitely correlated representations of $\E_n$ (or  of
$\O_n$ in the Cuntz case).

\begin{thm}\label{unitary_inv}
Let $A = (A_1,\dots,A_n)$ and $B=(B_1,\dots,B_n)$ be contractive
$n$-tuples on finite dimensional spaces $\V_A$ and $\V_B$ respectively.
Let $S=(S_1,\dots,S_n)$ and $T=(T_1,\dots,T_n)$ be their joint minimal
isometric dilations on Hilbert spaces $\H_A$ and $\H_B$;
and let $\sigma_A$ and $\sigma_B$ be the induced representations of
$\E_n$.
Let $\td{\V}_A$ be the subspace spanned by all minimal
$\fA^*$-invar\-iant subspaces $\W$ on which
$\sum_{i=1}^n A_iA_i^*|_\W = I_\W$; and similarly define $\td{\V}_B$.
Then $\sigma_A$ and $\sigma_B$ are unitarily equivalent if and only if
\begin{enumerate}
\item $\rank ( I_{\V_A} - \sum_{i=1}^n A_i A_i^* ) =
       \rank ( I_{\V_B} - \sum_{i=1}^n B_i B_i^* )$; \vspace{.2ex} and
\item $A^*|_{\td{\V}_A}$ is unitarily equivalent to $B^*|_{\td{\V}_B}$.
\end{enumerate}
\end{thm}

\Prf The two representations are equivalent if and only if they have
the same pure rank and the Cuntz parts are unitarily equivalent.  By
Theorem~\ref{general}, the algebra $\fS$ contains the projection onto
$\td{\V}_A$.  It is the unique maximal finite rank projection in $\fS$.
Therefore the restriction $A^*|_{\td{\V}_A}$ is a
unitary invariant.  Conversely, if these two conditions hold, then the
unitary identifying $A^*|_{\td{\V}_A}$ and $B^*|_{\td{\V}_B}$ extends
to a unitary equivalence between the dilations $S_A$ of
$\td{A} := P_{\td{\V}_A} A|_{\td{\V}_A}$ and $S_B$ of
$\td{B} :=  P_{\td{\V}_B} B|_{\td{\V}_B}$ because of the uniqueness of
the minimal isometric dilation.  This identifies the restriction of
$S_A$ to $\fS[\td{\V}_A] = \H_{Ac}$, namely the Cuntz part of $S_A$,
with the corresponding Cuntz part of $S_B$.  The pure rank condition
allows a unitary equivalence between the two pure parts.
\bx

Bratteli and Jorgensen \cite{BJendII} give a detailed analysis of
representations of the Cuntz algebra which has a lot in common with
our results.  They look somewhat different since they concentrate on
the state and not on the restriction to the subspace $\V$.  In particular,
their contractions are not the same as ours.  They point out the
relationship in the discussion preceding their Theorem~5.3.  They obtain
our Corollary~\ref{irreducible} in the Cuntz case, and in particular
recognize the role of the completely positive map $\Phi$.  Again
however, their different normalization results in a different map.
But they do not appear to classify these representations up to unitary
equivalence.  The reason they do not succeed is that they did not
identify the subspace which we call $\td{\V}$, and instead work with a
subspace they call $\V_k$ which is often strictly larger.  The space
$\td{\V}$ does not occur in their hierarchy of invariant subspaces.
Instead, they specialize in section~7 to a smaller class which they call
diagonalizable shifts.  These they do completely classify up to
unitary equivalence.  We have not determined in this case how their
special invariants relate to ours.

\begin{cor}\label{hyperreflexive}
The algebra $\fS$ determined by the joint isometric dilation of a
contractive $n$-tuple on a finite dimensional space is hyper-reflexive
with distance constant at most 5.
\end{cor}

\Prf This follows immediately from \cite[Theorem~3.14]{DP1} since the
algebra $\fS$ is unitarily equivalent to the algebra of certain atomic
representations.  Indeed, the projection $P = P_{\td{\V}}$ belongs to
$\fS$ and $\fS P = \fW P$ where $\fW$ is a type I von Neumann algebra
containing the projection $P$.  Thus by Christensen's result
\cite{Chr} which shows that type I von Neumann algebras have distance
constant at most 4, we obtain the same for our slice.  The upper bound
for the distance constant of $\fL_n$ was improved by Bercovici
\cite{Berc} to 3 from the original 51.  Arguing as in
\cite{DP1}, we obtain a distance constant no larger than
$(3^2+4^2)^{1/2} = 5$.
\bx

\section{Similarity}\label{S:similarity}

Now consider the question of when two contractive $n$-tuples are
similar, and the effect on their dilations.  The first step is to show
that the Cuntz parts must be unitarily equivalent.  Thus the question
of similarity reduces to the pure parts.  First we need a variant of
Lemma~\ref{expn1}.

\begin{lem}\label{expn3}
Suppose that an $n$-tuple $(A_1,\dots,A_n)$ acts on a finite dimensional space
$\V$, and generates $\B(\V)$ as an algebra.  Moreover suppose that
$\Phi(X) = \sum_{i=1}^n A_i X A_i^*$ is unital.  Then the only
self-adjoint operators $X$ satisfying $\Phi(X)\le X$ are scalar, and in
particular are fixed points.
\end{lem}

\Prf Since $\Phi(I)=I$, we may translate $X$ so that $X\ge 0$ and $0$
belongs to its spectrum.  Let $\M=\ker X$.  This is a non-zero subspace.
Let $x \in \M$.  Then
\begin{align*}
  0 &= (\Phi(0)x,x) \le (\Phi(X)x,x)\\
    &= \sum_{i=1}^n (A_i X A_i^*x,x) = \sum_{i=1}^n \|X A_i^*x\|^2
     \le (Xx,x) = 0 .
\end{align*}
It follows that $\M$ is invariant for each $A_i^*$.  But by hypothesis,
the $A_i^*$'s generate the full matrix algebra, and thus have no
proper invariant subspaces.  So $\M=\V$ and $X=0$ is scalar.
\bx

\begin{cor}\label{similar_Cuntz}
Suppose that $A = (A_1,\dots,A_n)$ and $B = (B_1,\dots,B_n)$ are similar
contractive $n$-tuples in the finite dimensional algebra $\B(\V)$.
Let $\td{\V}_A$ and $\td{\V}_B$ denote the subspaces spanned by the
minimal $\fA^*$ and $\fB^*$-invar\-iant subspaces $\M$ on which
$AA^*|_\M = I_\M$ and $BB^*|_\M = I_\M$, respectively.
Then $P_{\td{\V}_A}A|_{\td{\V}_A}$ and $P_{\td{\V}_B}B|_{\td{\V}_B}$ are
unitarily equivalent.
\end{cor}

\Prf Let $T$ be the similarity such that $B=TAT^{-1}$.  Then
$B^* = T^{*-1}A^*T^*$.  So $T^{*-1}$ carries $\fA^*$-invar\-iant
subspaces onto $\fB^*$-invar\-iant subspaces.  Also $T^{*-1}$ preserves
minimality.  However it is not immediately evident that it preserves
the condition that $AA^*|_{\td{\V}_A} = I_{\td{\V}_A}$.

Let $\M$ be a minimal $\fA^*$-invar\-iant subspace of $\td{\V}_A$
on which $AA^*|_\M = I_\M$.
Then $T^{*-1}\M = \N$ is invariant for $\fB^*$.  It follows that if
$\bar{A}^*_i$,and $\bar{T}^{*-1}$ are the restrictions of $A_i^*$ and
$T^{*-1}$ to $\M$, and $\bar{B}^*_i$ is the restriction of $B_i^*$ to
$\N$, then $\bar{B}_i^* = \bar{T}^{*-1} \bar{A}^*_i \bar{T}^*$.  Let
$\bar{\Phi}(X) = \sum_{i=1}^n \bar{A}_i X \bar{A}_i^*$ on $\B(\M)$.
It is easy to verify that $\bar{\Phi}$ is unital.

Now compute that
\[
  I_\N \ge \sum_{i=1}^n \bar{B}_i \bar{B}_i^* = \sum_{i=1}^n
     \bar{T} \bar{A}_i \bar{T}^{-1} \bar{T}^{*-1} \bar{A}_i^* \bar{T}^*
  = \bar{T} \Phi( \bar{T}^{-1} \bar{T}^{*-1} ) \bar{T}^* .
\]
Therefore
$\Phi( \bar{T}^{-1} \bar{T}^{*-1} ) \le \bar{T}^{-1} \bar{T}^{*-1}$.
By Lemma~\ref{expn3}, it follows that $\bar{T}^{-1} \bar{T}^{*-1}$ is
scalar.  So up to a scaling factor, $\bar{T}$ is unitary.

This shows that the restrictions of $\fA^*$ to each minimal
$\fA^*$-invar\-iant subspace $\M$ of $\td{\V}_A$ on which
$AA^*|_\M = I_\M$ is unitarily
equivalent to the corresponding subspace of $\fB^*$.  Since $\td{\V}_A$
and $\td{\V}_B$ are each the orthogonal direct sum of such subspaces, it
follows that the restriction to these larger subspaces are unitarily
equivalent (although $T$ itself need not be a multiple of a unitary on
the whole space).  Thus $A^*|_{\td{\V}_A}$ is unitarily equivalent to
$B^*|_{\td{\V}_B}$.  Equivalently, the compressions
$P_{\td{\V}_A}A|_{\td{\V}_A}$ and $P_{\td{\V}_B}B|_{\td{\V}_B}$
are unitarily equivalent.
\bx

\begin{cor}\label{Cuntz_sim}
Suppose that $A = (A_1,\dots,A_n)$ and $B = (B_1,\dots,B_n)$ are similar
contractive $n$-tuples in the finite dimensional algebra $\B(\V)$.
Let $S_i$ and $T_i$ be their respective minimal joint isometric
dilations.  Then the Cuntz parts of $S_i$ and $T_i$ are unitarily
equivalent.
\end{cor}

\Prf This is immediate from the proposition above and the fact that the
Cuntz part of $S_i$ and $T_i$ are determined by the compressions of
$A_i$ and $B_i$ to the subspaces $\td{\V}_A$ and $\td{\V}_B$
respectively by Corollary~\ref{unitary_inv}.
\bx

\begin{eg}
Now we show through a couple of examples that the pure part of the
dilation is not preserved by similarity.  This first example shows that
one dilation can be strictly Cuntz type while a similarity can
introduce a pure part.  Consider
\[
  A_1 = \begin{bmatrix}1&0\\0&0\end{bmatrix} \qqand
  A_2 = \begin{bmatrix}0&0\\1&0\end{bmatrix} .
\]
This is of Cuntz type since $A_1A_1^*+A_2A_2^* = I$.
Moreover there is a unique minimal $\fA^*$-invar\-iant subspace,
$\bbC e_1$.  The dilation of this pair is thus irreducible by
Corollary~\ref{irreducible}, and is determined by the 1-dimensional
restrictions 1 and 0 of $A_1^*$ and $A_2^*$ to $\bbC e_1$.
In fact, this is easily seen to be the atomic representation
$\sigma_{1,1}$.

However, this pair is similar via
$T = \begin{bmatrix}1&0\\0&1/2\end{bmatrix}$ to
\[
  B_1 = \begin{bmatrix}1&0\\0&0\end{bmatrix} \qqand
  B_2 = \begin{bmatrix}0&0\\1/2&0\end{bmatrix} .
\]
The restrictions of $B_i^*$ to the unique minimal $\fB^*$-invar\-iant
subspace are still 1 and 0 respectively; and they determine a
dilation which has the representation $\sigma_{1,1}$ as a summand.
However, since
\[
  \rank(I - B_1 B_1^* - B_2 B_2^*) =
  \rank \begin{bmatrix}0&0\\0&3/4\end{bmatrix} = 1 ,
\]
the pure rank of this representation is 1.
\end{eg}

\begin{eg}
A second easy example shows that even in the pure case, the pure rank
is not a similarity invariant.  Fix an orthonormal basis
$e_1,\dots,e_n$ for $\V$.
Let $A_1 = \frac12 e_1e_1^*$ and $A_i = e_ie_1^*$ for $2 \le i \le n$.
Then $I - \sum_{i=1}^n A_iA_i^* = \frac34 e_1e_1^*$ is rank 1, and its
range $\bbC e_1$ is $\fA$-cyclic.  So by Corollary~\ref{irreducible},
this yields an irreducible pure dilation.

However, this is similar via $T = I + e_1e_1^*$ to
$B_1 = A_1$ and $B_i = \frac12 A_i$ for $2 \le i \le n$.  This
$n$-tuple satisfies $I - \sum_{i=1}^n B_iB_i^* = \frac34 I$, which has
rank $n$.  So this dilation has pure rank $n$.
\end{eg}

We wish to provide more detail about the effect of similarity on pure
representations.  By Popescu\cite{Pop_diln}, the dilation is
pure if and only if
\[
  \wotlim\limits_{k\to\infty} \sum\limits_{|w|=k} A_w A_w^* = 0 .
\]
He calls these $n$-tuples $C_0$ contractions, and provides a
\wot-continuous functional calculus in \cite{Pop_vN}.
We can analyze this using the theory of representations of
$\fL_n$ developed in \cite{DP2}.

Let $A=(A_1,\dots,A_n)$ be a $C_0$-contraction on $\V$ with pure
minimal isometric dilation $S_i \simeq L_i^{(s)}$.  This determines a
\wot-continuous representation $\Phi_A$ of $\fL_n$ which sends $X$ to
$P_\V X^{(s)}|_\V$.
In particular, $\Phi_A(L_w) = A_w := A_{i_1}\dots A_{i_k}$ for every
word $w = i_1 \dots i_k$ in $\F_n$.
The kernel $\fJ = \ker \Phi_A$ is a \wot-closed ideal of $\fL_n$.
By \cite[Theorem 2.1]{DP2}, this ideal is determined by its range
$\M = \ol{\fJ \K_n}$, which is an invariant subspace for both $\fL_n$
and its commutant $\fR_n$.
The representation of compression of $\fL_n$ to $\M^\perp$ has the same
kernel.  We wish to determine to what extent $A$ can be recovered from
the compression of $L$ to $\M^\perp$.

To get a feeling for the situation, consider the case in which the $A_i$
are $d\times d$ matrices which generate $\fM_d$ as an algebra.  Then
$\Phi_A$ maps $\fL_n$ onto $\fM_d$.  The kernel $\fJ$ will then have
codimension $d^2$, and therefore the dimension of $\M^\perp$ is also
$d^2$. The compression homomorphism to $\M^\perp$ factors through
$\Phi_A$. Since $\fM_d$ has only one irreducible representation up to
similarity, the compression to $\M^\perp$ must be similar to the direct
sum of $d$ copies of $\Phi_A$.  In particular, $\M^\perp$ will
decompose into a (non-orthogonal) direct sum of $d$ subspaces which are
$\fL_n^*$-invar\-iant such that the compression of $L$ is similar to $A$.

Nevertheless, $\Phi_A$ need not occur as a compression of $L$ to some
$\fL_n^*$-invar\-iant subspace.  This could occur only if $\Phi_A$ has
pure rank 1, which need not be the case.  However, this shows that there
are representations similar to $\Phi_A$ which do have pure rank 1.
Moreover it turns out that in a certain sense, these similarities of
pure rank 1 are the extreme points of those representations similar to
$\Phi_A$.  This will be established by showing that $\Phi_A$ can be
recovered as a C*-convex combination of pure rank 1 representations.

\begin{thm}\label{C_0_reps}
Let $A=(A_1,\dots,A_n)$ be a $C_0$-contraction on a $d$-dimensional
space $\V$.
Let $\fJ$ be the kernel of the \wot-continuous representation
$\Phi_A$ of $\fL_n$ that it determines.
Then $\Phi_A$ is unitarily equivalent to the compression of $\fL_n$ to
a semi-invar\-iant subspace $\S=\N_1\ominus\N_2$, where $\N_1=\fL_n[\S]$
and $\N_2 = \N_1\ominus\S$ belong to $\Lat\fL_n$.\vspace{.2ex}

Let $\M = \ol{\fJ\K_n}$ be the corresponding
$\fL_n$ and $\fR_n$-invar\-iant subspace associated to $\fJ$.
Then there are at most $d$ wandering vectors $\zeta_j$, say for
$1\le j\le s$ where $s\le d$, with $\fL_n[\zeta_i]$ pairwise orthogonal,
such that
\[
 \N_1 = \sum_{j=1}^s\upplus R_{\zeta_j}\K_n \qqand
 \N_2 \supset \sum_{j=1}^s\upplus R_{\zeta_j}\M .
\]

Moreover, the subspaces $\M_j = R_{\zeta_j}^* \S$ are
$\fL_n^*$-invar\-iant subspaces of $\M^\perp$, and $\dim(\M_j)\le d$.
The contractive $n$-tuples $B_j = (B_{j1},\dots,B_{jn})$ obtained by
compression of $L$ to $\M_j$ have pure rank 1 and
$\ker \Phi_{B_j} \supset \fJ$.
There is an isometry $X$ mapping $\S$ into $\sum^{\oplus s}_{j=1} \M_j$
so that
\[
  \Bigl( \sum^s_{j=1}\upplus B_{ji}^* \Bigr) X = X A_i^*
  \qfor 1 \le i \le n .
\]
Thus $A^*$ is unitarily equivalent to the restriction
of $\sum^{\oplus s}_{j=1} B_j^*$ to an invariant subspace.
Consequently, $A = X^* \sum^{\oplus s}_{j=1} B_j X$ is a C*-convex
combination of the $B_j$'s.

In particular, when $\fA = \B(\V)$, each subspace $\M_j$ is
$d$-dimensional and each $n$-tuple $B_j$ is similar to $A$.
\end{thm}

\Prf  The isometric dilation $S=(S_1,\dots,S_n)$ of $A$ has pure rank
$s = \rank (I - \sum_{i=1}^n A_i A_i^*) \le d$.
The identity representation of $\fL_n$ contains many invariant
subspaces with infinite dimensional wandering space; and thus an
infinite multiple of the identity representation is contained in
$\fL_n$.  So we may assume that $\Phi_A$ is the compression of $\fL_n$
to a semi-invar\-iant subspace $\S$ of $\fL_n$ itself.
The minimal choice of a pair of $\fL_n$-invar\-iant subspaces with
difference $\S$ is given by $\N_1=\fL_n[\S]$ and $\N_2 = \N_1\ominus\S$
\cite{Sar}.

Now $\N_1$ has a wandering space $\W$ of dimension $s$.
Choose an orthonormal basis $\zeta_j$, $1\le j\le s$, for $\W$.
By \cite[Theorem 2.1]{DP1}, there is an isometry $R_{\zeta_j}$ in
$\fR_n$ with range equal to the cyclic $\fL_n$-invar\-iant subspace
$\fL_n[\zeta_j]$.  Then $\N_1 = \sum_{j=1}^{\oplus s} R_{\zeta_j}\K_n$.
Since the kernel of the compression to $\S$ is $\fJ$,
\begin{align*}
 \N_2 &\supset \fJ \S = \fJ \fL_n \S = \fJ \N_1 \\
  &= \sum_{j=1}^s \fJ R_{\zeta_j} \K_n
   = \sum_{j=1}^s R_{\zeta_j} \fJ \K_n = \sum_{j=1}^s R_{\zeta_j} \M .
\end{align*}

The subspaces $\M_j = R_{\zeta_j}^* \S$ are contained in
$\M^\perp$, and have dimension at most $d = \dim \S$.
Moreover, they are $\fL_n^*$-invar\-iant because of the identity
\[
 \fL_n^* R_{\zeta_j}^* \S = R_{\zeta_j}^* \fL_n^* \S \subset
 R_{\zeta_j}^* \N_2^\perp = R_{\zeta_j}^*\S .
\]
Let $B_j$ denote the contractive $n$-tuple obtained by compression of $L$
to $\M_j$.  Clearly $L$ is an isometric dilation of $B_j$.  The minimal
dilation is obtained by restricting $L$ to $\fL_n[\M_j]$.
However by Lemma~\ref{reducing}, this is a reducing subspace of $\K_n$.
Since the commutant $\fR_n$ of $\fL_n$ contains no idempotents
\cite[Corollary~1.8]{DP1}, this space must be all of $\K_n$.
Thus the $n$-tuple $B_j$ has pure rank 1 for each $1 \le j \le s$.
Since $\M_j$ is contained in $\M^\perp$, it follows that $\ker \Phi_{B_j}$
contains the ideal $\fJ$.

Now notice that $R_{\zeta_j} \M_j$ are pairwise orthogonal subspaces, and
\[
 \sum_{j=1}^s\upplus R_{\zeta_j} \M_j =
 \sum_{j=1}^s\upplus R_{\zeta_j} R_{\zeta_j}^* \S \supset
 \Bigl( \sum_{j=1}^s R_{\zeta_j} R_{\zeta_j}^* \Bigr) \S =
 P_{\N_1} \S = \S .
\]
This allows us to identify $\S$ isometrically with a subspace of
$\sum_{j=1}^{\oplus s} \M_j$.
Let $X_j = R_{\zeta_j}^* P_\S$ be considered as a map from $\S$ into $\M_j$.
Define $X$ to be the column matrix
$\begin{bmatrix} X_1& \dots& X_s \end{bmatrix}^t$.  Then
\[
 X^*X = \sum_{j=1}^s P_\S R_{\zeta_j} R_{\zeta_j}^* P_\S
      = P_\S P_{\N_1}P_\S = P_\S .
\]
So $X$ is an isometry of $\S$ into $\sum_{j=1}^{\oplus s} \M_j$.
One may compute
\[
  R_{\zeta_j}^* L_i^* P_\S =  \bigl( R_{\zeta_j}^* P_{\N_1} \bigr)
                              \bigl( P_{\N_2}^\perp L_i^* P_\S \bigr) =
 R_{\zeta_j}^* P_\S L_i^* P_\S .
\]
Therefore identifying $A_i$ with $P_\S L_i P_\S$, we obtain
\begin{align*}
  \bigl( \sum_{j=1}^s\upplus B_{ji}^* \bigr) X &=
  \sum_{j=1}^s L_i^* R_{\zeta_j}^* P_\S =
  \sum_{j=1}^s R_{\zeta_j}^* L_i^* P_\S  \\ &=
  \sum_{j=1}^s R_{\zeta_j}^* P_\S L_i^* P_\S = X A_i^* .
\end{align*}
>From this it is evident that the range of $X$ is invariant for
$\sum_{j=1}^{\oplus s} B_{ji}^*$, and implements a unitary equivalence
between $A_i^*$ and this restriction of $\sum_{j=1}^{\oplus s} B_{ji}^*$.
Consequently, $X^* \bigl( \sum_{j=1}^{\oplus s} B_{ji} X \bigr) = A_i$
for $1 \le i \le n$.  This expresses $A$ as a C*-convex combination of
the pure rank one contractions $B_j$.

When $\fA$ is isomorphic to $\B(\V)\simeq \fM_d$,
then $\fL_n/\fJ$ is likewise isomorphic to $\fM_d$ and the compression of
$\fL_n$ to $\M^\perp$ is a representation of $\fM_d$ on a subspace of
dimension $d^2$.  The only representations of $\fM_d$ are
multiples of the identity representation up to similarity,
and the compression to $\M^\perp$ has multiplicity $d$.
Thus the $\fL_n^*$-invar\-iant subspaces of $\M^\perp$ have dimension
which is a multiple of $d$.
As $\M_j$ are non-zero and have dimension at most $d$, they are all
exactly $d$-dimensional and each $B_j^*$ is similar to $A^*$, whence
$B_j$ is similar to $A$.
\bx

Say that the $n$-tuple $A$ is {\it irreducible} if it generates
$\fA = \B(\V)$, or equivalently $\fA$ has no proper invariant subspaces.
When $A$ is an irreducible $C_0$-contraction, we see that the
compression representation to $\M^\perp$ takes the generators to an
$n$-tuple similar to the direct sum of $d$ copies of $A$.  In
particular, this occurs if $\|A\| < 1$.  So we obtain a complete
similarity invariant for an {\it arbitrary} irreducible $n$-tuple of
matrices (after scaling appropriately).

Restricting to the irreducible case is not just a matter of
convenience.  Simple examples show that multiplicity cannot be
detected from the set of polynomial identities that an $n$-tuple
satisfies.  For example, with $n=1$, take $A = J_2 \oplus 0^{(3)}$ and
$B = J_2^{(2)} \oplus 0$ where $J_2$ is the $2 \times 2$ nilpotent
Jordan matrix and 0 is a one-dimensional zero.  These two matrices
satisfy exactly the same polynomial identities.  The natural way to
distinguish them is to use rank.  Indeed, familiar invariants for
similarity of single matrices shows that the ranks of various
polynomials can be used to determine the multiplicity function.

\begin{cor}\label{sim_inv}
Suppose that $A=(A_1,\dots,A_n)$ and $B=(B_1,\dots,B_n)$ are two
$n$-tuples of $d\times d$ matrices which are irreducible
and strictly contractive, $\|A\|<1$ and $\|B\|<1$.  Then $A$ and $B$
are similar if and only if $\ker \Phi_A = \ker \Phi_B$.
\end{cor}

\Prf Clearly two similar $n$-tuples give rise to representations with
the same kernel.  Conversely, if they are irreducible, the kernel
determines the subspace $\M$.  We adopt the notation from the proof of
Theorem~\ref{C_0_reps}.  A minimal $\fL_n^*$-invariant subspace
$\M_j$ of $\M^\perp$ yields a compression representation $\Phi_j$ which
is similar to $\Phi_A$ by Theorem~\ref{C_0_reps}.
Likewise, $B$ determines the same subspaces, and thus $\Phi_B$ is also
similar to $\Phi_j$, and hence to $\Phi_A$.
\bx

The similarity question for $n$-tuples of matrices is an old one, and
the solution is complicated.
Friedland \cite{Fried} provides an algorithm for checking whether
two $n$-tuples $A=(A_1,\dots,A_n)$ and $B=(B_1,\dots,B_n)$ of $d\times
d$ matrices  are similar.  This is quite involved even for two
$2\times 2$ matrices, which he calculates explicitly.  The situation
simplifies when the two matrices are not simultaneously
triangularizable---which in the $2\times 2$ case is the same as
irreducibility.  In this case, the pairs $A$ and $B$ are similar if and
only if these five identities hold:
\begin{alignat*}{2}
  \Tr(A_1) &= \Tr(B_1) &\qquad \Tr(A_1^2) &= \Tr(B_1^2)\\
  \Tr(A_2) &= \Tr(B_2) &\qquad \Tr(A_2^2) &= \Tr(B_2^2)\\
  \Tr(A_1A_2) &= \Tr(B_1B_2) .&&
\end{alignat*}
In general, there is no explicit list of polynomials to check.

In our case, we do obtain a fairly small {\it finite} list of
invariants for an irreducible $n$-tuple.  Unfortunately, at this
point, we do not have an explicit method for computing these
invariants.  Nor are they polynomials.  The natural invariants in our
setting are isometries in $\fL_n$.  Polynomials can be obtained by a
simple approximation argument, but are no longer canonical.  In
the case of two $2\times 2$ matrices, we obtain exactly five
conditions.

\begin{thm}\label{sim_invariants}
Let $A=(A_1,\dots,A_n)$ be an irreducible $n$-tuple of $d\times d$
matrices with $\|A\|<1$.  The ideal $\fJ = \ker \Phi_A$ determines its
range space $\M = \ol{\fJ \K_n}$ with wandering dimension $1+(n-1)d^2$.
Thus there are $1+(n-1)d^2$ isometries $X_j$ in $\fL_n$ so that an
$n$-tuple $B$ of $d\times d$ matrices with $\|B\|<1$ is similar to $A$
if and only if $\Phi_B(X_j)=0$ for $1 \le j \le 1+(n-1)d^2$.

Moreover there is a set of $m = 1+(n-1)d^2$ polynomials $p_j$
in $n$ non-commuting variables such that an $n$-tuple $B$ of
$d\times d$ matrices with $\|B\|<1$ is similar to $A$ if and only if
$p_j(B)=0$ for $1 \le j \le m$.
\end{thm}

\Prf The space $\M$ has the same codimension as $\fJ$, which is $d^2$
since $\fL_n/\fJ$ is isomorphic to $\Alg\{A_1,\dots,A_n\} = \fM_d$.
Its wandering space is
{\allowdisplaybreaks
\begin{align*}
  \W &= \M \ominus \sum_{i=1}^n L_i \M \\
     &= \K_n \ominus \bigl( \M^\perp \oplus
                \sum_{i=1}^n L_i (\K_n \ominus \M^\perp) \bigr) \\
     &= \Bigl( \bigl( \K_n \ominus \sum_{i=1}^n L_i \K_n \bigr)
               \oplus \sum_{i=1}^n L_i \M^\perp \Bigr)
        \ominus \M^\perp \\
     &= \Bigl( \bbC \xi_e \oplus \sum_{i=1}^n L_i \M^\perp \Bigr)
        \ominus \M^\perp .
\end{align*}
}
This has dimension $m = 1 + (n-1) \dim \M^\perp = 1 + (n-1)d^2$.

Now $\M$ is invariant for both $\fL_n$ and its commutant $\fR_n$.
Since it is the latter, it decomposes \cite[Theorem~2.1]{DP1} as the
direct sum of $m$ cyclic $\fR_n$-invariant subspaces; and each is
the range of an isometry $X_j$ in $\fL_n$.  Thus by
\cite[Lemma~2.5]{DP2}, we obtain that $\fJ = \sum_{j=1}^m X_j \fL_n$.

Therefore $\ker \Phi_B$ contains $\fJ$ if and only if $\Phi_B(X_j) = 0$
for $1 \le j \le m$.  Moreover, since $A$ is irreducible, $\fJ$ is a
maximal ideal.  Thus this condition ensures that $\ker\Phi_B = \fJ$.
In particular, $\Alg\{B_1,\dots,B_n\}$ is isomorphic to $\fL_n/\fJ
\simeq \fM_d$; and hence $B$ is also irreducible.  Therefore $B$ and
$A$ are similar by Corollary~\ref{sim_inv}.

To obtain polynomials, we notice that the algebra $\fP$ which is the
algebraic span of $\{L_w : w \in \F_n \}$ is \wot-dense in $\fL_n$.
Let $\fI = \fJ \cap \fP$ be the ideal of all polynomials which
annihilate $A$.  The algebra (without closure) generated by the
$A_i$'s is $\fM_d$.  So the map $\Phi_A$ takes $\fP$ onto $\fM_d$ with
kernel $\fI$; and takes $\fL_n$ onto $\fM_d$ with kernel $\fJ$.  It
follows from the Hahn-Banach theorem that $\fI$ is \wot-dense in $\fJ$.

Let $\eps = (1+nd^2)^{-1}$.  For each $1 \le j \le m$, choose
polynomials $p_j\in \fI$ such that $\| p_j(L) - X_j \| < \eps$.
We claim that $p_j(L)$ generate $\fJ$ as a norm-closed right ideal.
For let $J\in\fJ$.  By \cite[Lemma~2.5]{DP2}, there are elements $Y_j
\in \fL_n$ such that $J = \sum_{j=1}^{m} X_j Y_j =: XY$.
Moreover, the row operator
$X = \bigl[ X_1 \; \dots \; X_m \bigr]$ is an isometry.
Hence the column operator
$Y = \bigl[ Y_1 \; \dots \; Y_m \bigr]^t$
has $\|Y\| = \|J\|$.
Let $P = \bigl[ p_1(L) \; \dots \; p_m(L) \bigr]$.
It follows that
\[
  \| J - PY \| \le \| X - P \| \|Y\| < \frac{m}{1+nd^2} \|J\| .
\]
Since $m/(1+nd^2)<1$, the right ideal generated by $P$
is norm dense in $\fJ$ as claimed.

Therefore the condition that $p_j(B)=0$ for $1 \le j \le m$
is equivalent to the condition $\Phi_B(X_j) = 0$, and thus is
equivalent to joint similarity to $A$.
\bx

While the $X_j$'s are needed to generate $\fJ$ as a \wot-closed {\it
right} ideal, there will generally be redundancies as generators for
$\fJ$ as a two-sided ideal.
So $1 + (n-1)d^2$ is an upper bound on the number of
test elements needed.  It would be interesting to have better bounds
on the number of generators for a two-sided ideal.

We observe that the existence of a determining set of polynomials for
an irreducible $n$-tuple can be deduced directly by elementary means.
One can write down polynomials in $A$ representing the matrix units of
$d\times d$ matrices and their relations.  In fact $O(d^2)$ generators
and relations suffice.  Then each $A_i$ can be expressed as a
combination of matrix units.  This requires only $n + O(d^2)$
polynomials, which is somewhat better than our bound.  In many
concrete cases, this simple bare hands approach is the best.

On the other hand, our result provides an algorithm for obtaining a
set of generators for the ideal $\fJ$.  Perhaps this will prove to be
of some use.

\begin{eg}
This example illustrates parts of the previous two theorems.
Consider the pair of $2 \times 2$ matrices
\[ A_1 =
 \begin{bmatrix} 0& 1/2\\ 1/2 & 0 \end{bmatrix}
     \qqand
 A_2 = \begin{bmatrix}1/2&0\\0&0\end{bmatrix} .
\]
Since $I - A_1A_1^* - A_2A_2^* =
  \begin{sbmatrix}1/2&0\\0&3/4\end{sbmatrix}$ has \vspace{.2ex}
rank 2, this determines a pure isometric dilation of pure rank 2.
The algebra $\fA = \fM_2$, and thus the representation $\Phi_A$ of
$\fL_2$ is irreducible.  The kernel will be a \wot-closed maximal
ideal $\fJ$ of codimension 4.
Therefore the subspace $\M = \ol{ \fJ \K_2}$ will be codimension 4.

The matrices $A_1$ and $A_2$ satisfy certain relations that express the
fact that
\[
  \fM_2 = \Alg \{ A_1, A_2 \} = \spn \{ I, A_1, A_2, A_1A_2 \} .
\]
A natural and sufficient list is
{\renewcommand{\labelenumi}{(\arabic{enumi})}
\begin{enumerate}
 \item    \qquad $4A_1^2 = I$
 \item    \qquad $2A_2^2 = A_2$
 \item    \qquad $A_2A_1A_2 = 0$
 \item    \qquad $8A_1A_2A_1 = I-2A_2$
 \item    \qquad $2A_1A_2 + 2A_2A_1 = A_1$
\end{enumerate}}
\noindent However, (5) and (2) can be derived from the others.
So (1), (3) and (4) are sufficient.

The ideal $\fJ$ is therefore generated as a two-sided ideal by the set
\[ \J = \{ I-4L_1^2,\; L_2L_1L_2,\; I-2L_2-8L_1L_2L_1 \} ,\]
since the quotient will be $\fM_2$.  Therefore the range $\M$ is the
$\fL_2 \fR_2$-invar\-iant subspace generated by $\J \xi_e$,
\begin{align*}
  \M &= \ol{ \fJ \K_2} = \ol{ \fL_2 \fR_2 \J \xi_e} \\ &=
  \spn \{ \xi_{uv}-4\xi_{u11v}, \; \xi_{u212v}, \;
           \xi_{uv}-2\xi_{u2v}-8\xi_{u121v} : u,v\in\F_2 \} .
\end{align*}

We wish to determine $\M^\perp$.  To this end, define $\Omega_1$ to be
the set of all words in $11=1^2$ and $2$,
\[ \Omega_1 = \{ w = 1^{2k_0}21^{2k_1}2\dots 21^{2k_s} :
     k_i\ge0,s\ge0 \} .\]
Let $\Omega_2 = \Omega_1 1$, $\Omega_3 = 1 \Omega_1$ and
$\Omega_4 = \{e\} \cup 1 \Omega_1 1$.  Define
\[
  x_i = \sum_{w\in\Omega_i} 2^{-|w|} \xi_w   \qfor 1\le i\le 4.
\]
Then a computation shows that $\M^\perp = \spn\{x_1,x_2,x_3,x_4\}$.
These vectors are not orthogonal, nor of constant length.  Indeed,
\[
 \|x_1\|^2 = 16/11\, , \qquad
 \|x_2\|^2 = \|x_3\|^2 = 4/11  \qand
 \|x_4\|^2 = 12/11 .
\]
The pair $\{x_1,x_4\}$ is orthogonal to $\{x_2,x_3\}$,
but
\[
 (x_1, x_4) = 16/15 \qqand (x_2, x_3) = 4/15 .
\]

Another matrix calculation relative to the ordered basis
$\{x_1, \ldots,x_4\} $ for $\M^\perp$ shows that
\[
  L_1^*|_{\M^\perp}=
   \begin{bmatrix}
    0&0&1/2&0\\
    0&0&0&1/2\\
    1/2&0&0&0\\
    0&1/2&0&0
   \end{bmatrix} \qand
  L_2^*|_{\M^\perp}=
   \begin{bmatrix}
    1/2&0&0&0\\
    0&1/2&0&0\\
    0&0&0&0\\
    0&0&0&0
   \end{bmatrix}.
\]
Let $\Phi_c$ denote the representation of compression to $\M^\perp$.
This calculation shows that $\Phi_c$ is similar (but not
unitarily equivalent) to the direct sum of two copies of $\Phi_A$.
Thus the compression to any two dimensional $\fL_2^*$-invar\-iant
subspace of $\M^\perp$ is similar to $\Phi_A$.
As noted in the proof of Theorem~\ref{C_0_reps}, these representations all
have pure rank 1.  In particular, $\Phi_A$ does not occur as such a
compression.
It is also a fact that $\Phi_c$ is not unitarily equivalent to an
orthogonal direct sum of two representations.

Next, we compute the 2-dimensional $\fL_2^*$-invar\-iant subspaces
of $\M^\perp$.   Examining the representation $\Phi_c$, we observe
that the two subspaces $\M_1 := \spn \{ x_1, x_3 \}$ and
$\M_2 := \spn \{ x_2, x_4 \}$ are invariant for $\fL_2^*$.
Setting  $\eta_i = x_i/ \|x_i\|$, we find that $\{\eta_1,\eta_3\}$ and
$\{\eta_2,\eta_4\}$ are orthonormal bases for $\M_1$ and $\M_2$
respectively.  Compute
\begin{alignat*}{2}
 &L_1^*|_{\M_1} \simeq B_1^* =
 \begin{bmatrix} 0&1\\ 1/4 &0 \end{bmatrix} &\qqand
 &L_2^*|_{\M_1} \simeq B_2^* =
 \begin{bmatrix} 1/2 &0\\ 0&0 \end{bmatrix} ,\\
 &L_1^*|_{\M_2} \simeq C_1^* =  \begin{bmatrix}
 0&1/\sqrt{12} \\ \sqrt{3}/2 &0 \end{bmatrix}
 &\qqand &L_2^*|_{\M_2} \simeq C_2^* =
 \begin{bmatrix} 1/2 &0\\ 0&0 \end{bmatrix} .
\end{alignat*}
These must be pairs which have pure rank 1, as is verified by
computing the ranks of
\[
  I \!-\! B_1B_1^* \!-\! B_2B_2^* =
   \begin{bmatrix} 11/16&0\\ 0&0 \end{bmatrix} \qand
  I \!-\! C_1C_1^* \!-\! C_2C_2^* =
   \begin{bmatrix} 0&0\\ 0&11/12 \end{bmatrix} .
\]

The representation $\Phi_c$ factors through a representation of $\fM_2$
of multiplicity 2.
Thus every 2-dimensional $\fL_2^*$-invar\-iant subspace is the
cyclic subspace determined by its intersection with the range of
$L_2^*|_{\M^\perp} = \spn \{ \eta_1, \eta_2 \}$, namely
$\bbC (\alpha \eta_1 + \beta \eta_2)$ for
$|\alpha|^2 + |\beta|^2 = 1$.  The second vector spanning the
subspace must be the image under $2L_1^*$,
\[
  2L_1^* (\alpha \eta_1 + \beta \eta_2) =
  2\alpha B_1^* \eta_1 + 2\beta C_1^* \eta_2 =
  \tfrac{\alpha}{2} \eta_3 + \sqrt{3} \beta \eta_4 .
\]
 A typical subspace of this form is
\[
  \M_{\alpha,\beta} = \spn \{ \alpha \eta_1 + \beta \eta_2,
   \tfrac{\alpha}{2} \eta_3 + \sqrt{3} \beta \eta_4 \} .
\]
However it is sufficient just to use $\M_1$ and $\M_2$, as they
correspond to a particular choice of a basis for the wandering
subspace of $\fL_2[\S]$.

Since rank considerations show that the subspace $\S$ cannot be
$\fL_2^*$-invar\-iant, we must write $\S$ as the difference of two
$\fL_2$-invar\-iant subspaces of multiplicity 2.  We will look for an
orthonormal set $\{\zeta_1,\zeta_2\}$ to span $\S$ of the form
\[
  \zeta_1 = \alpha R_1\eta_1 + \beta R_2\eta_2 \qand
  \zeta_2 = \tfrac{\alpha}{2}R_1\eta_3 + \sqrt{3}\beta R_2\eta_4 .
\]
They are always orthogonal, so the condition that they be norm one
requires that
\[
 |\alpha|^2 + |\beta|^2 = 1 = \tfrac{1}{4} |\alpha|^2 + 3|\beta|^2 .
\]
This has the solution $|\alpha|^2 = 8/11$ and $|\beta|^2 = 3/11$.
Therefore set
\[
 \zeta_1 =
  \sqrt{8/11}R_1\eta_1 + \sqrt{3/11} R_2\eta_2 \;\;\text{and}\;\;
 \zeta_2 =
  \sqrt{2/11}R_1\eta_3 + \sqrt{9/11}R_2\eta_4 .
\]

Another computation shows that
\begin{alignat*}{2}
 L_1^*\zeta_1 &= \tfrac{1}{2}\zeta_2 + \tfrac{1}{\sqrt{2}}\xi_e
 &\quad\qquad  L_2^*\zeta_1 &= \tfrac{1}{2}\zeta_1
 \\  L_1^*\zeta_2 &= \tfrac{1}{2}\zeta_1
 &\quad\qquad  L_2^*\zeta_2 &= 0 + \tfrac{1}{2}\xi_e
\end{alignat*}
Thus we see that $\S:=\spn\{\zeta_1,\zeta_2\}$ is a semi-invar\-iant
subspace $\S = \N_1 \ominus \N_2$ where $\N_1 = \{\xi_e\}^\perp =
\fL_2[\xi_1,\xi_2]$ and $\N_2 = \{ \xi_e, \zeta_1, \zeta_2 \}^\perp$.
Moreover these identities show that the compression of
$L_i$ to $\S$ is unitarily equivalent to $A_i$ for $i=1,2$.

Let us compute the operators $X_1$ and $X_2$ promised in
Theorem~\ref{C_0_reps}.
The projection onto $\S$ is given by
$P_\S = \zeta_1\zeta_1^* + \zeta_2\zeta_2^*$.  Then
\[
 X_1 = R_1^* P_\S = \sqrt{\tfrac{8}{11}}\eta_1\zeta_1^* +
                    \sqrt{\tfrac{2}{11}}\eta_3\zeta_2^*
\]
and
\[
 X_2 = R_2^* P_\S = \sqrt{\tfrac{3}{11}}\eta_2\zeta_1^* +
                    \sqrt{\tfrac{9}{11}}\eta_4\zeta_2^* .
\]
So recalling the matrix forms for $B_i$ and $C_i$, we obtain
\begin{align*}
 \sum_{i=1}^2 X_i^* L_1 X_i &=
  \begin{bmatrix} \sqrt{8/11} & 0\\ 0 & \sqrt{2/11} \end{bmatrix}
  \begin{bmatrix} 0 & 1/4 \\ 1 & 0 \end{bmatrix}
  \begin{bmatrix} \sqrt{8/11} & 0\\ 0 & \sqrt{2/11} \end{bmatrix}
 \\&\quad +
  \begin{bmatrix} \sqrt{3/11} & 0 \\ 0 & \sqrt{9/11} \end{bmatrix}
  \begin{bmatrix} 0 & \sqrt{3}/2 \\ 1/\sqrt{12} & 0 \end{bmatrix}
  \begin{bmatrix} \sqrt{3/11} & 0 \\ 0 & \sqrt{9/11} \end{bmatrix}
 \\&=
  \begin{bmatrix} 0 & 1/11 \\ 4/11 & 0 \end{bmatrix}  +
  \begin{bmatrix} 0 & 9/22 \\ 3/22 & 0 \end{bmatrix}  =
  \begin{bmatrix} 0 & 1/2  \\ 1/2  & 0 \end{bmatrix}
 = A_1
\end{align*}
and
\begin{align*}
 \sum_{i=1}^2 X_i^*L_2X_i &=
  \begin{bmatrix} \sqrt{8/11} & 0 \\ 0 & \sqrt{2/11} \end{bmatrix}
  \begin{bmatrix} 1/2 & 0 \\ 0 & 0 \end{bmatrix}
  \begin{bmatrix} \sqrt{8/11} & 0\\ 0 & \sqrt{2/11} \end{bmatrix}
 \\&\quad +
  \begin{bmatrix} \sqrt{3/11} & 0 \\ 0 & \sqrt{9/11} \end{bmatrix}
  \begin{bmatrix} 1/2 & 0 \\ 0 & 0 \end{bmatrix}
  \begin{bmatrix} \sqrt{3/11} & 0\\ 0 & \sqrt{9/11} \end{bmatrix}
 \\&=
  \begin{bmatrix} 4/11 & 0 \\ 0 & 0 \end{bmatrix}  +
  \begin{bmatrix} 3/22 & 0 \\ 0 & 0 \end{bmatrix}  =
  \begin{bmatrix} 1/2  & 0 \\ 0 & 0 \end{bmatrix}
 = A_2 .
\end{align*}
\bigbreak

The wandering space for $\M$ has dimension 5, and as in the proof
of Theorem~\ref{sim_invariants} is given by
\[
  \W = \spn \{ \xi_e, L_1 x_j, L_2 x_j : 1 \le j \le 4 \}
       \ominus \spn \{ x_j : 1 \le j \le 4 \} .
\]
We do not compute a basis for this space, as the result is not
particularly illuminating.  But such a basis corresponds to 5 isometries
$X_1,\dots,X_5$ in $\fL_n$ such that $\fJ = \sum_{j=1}^5 X_j \fL_n$.
Thus $\ker \Phi_B = \fJ$ if and only if $\Phi_B(X_j) = 0$ for
$1 \le j \le 5$.  But as we noted earlier in our remarks, finding
generators for $\fJ$ as a {\it right} ideal is overkill.  It suffices
to use generators for $\fJ$ as a two sided ideal.  Thus it is
sufficient to verify the 3 polynomial conditions:
{\renewcommand{\labelenumi}{(\arabic{enumi})}
\begin{enumerate}
 \item    \qquad $4 B_1^2 = I$
 \stepcounter{enumi}
 \item    \qquad $B_2 B_1 B_2 = 0$
 \item    \qquad $8 B_1 B_2 B_1 = I - 2 B_2$
\end{enumerate}}\noindent
It is an easy exercise to verify directly that these conditions
suffice to determine the pair up to similarity.

Pursuing this example further, let us consider other contractive pairs
which are similar to $A$.  A calculation shows that these
pairs are unitarily equivalent to a pair of the form
\[
 \td{A}_1 = \begin{bmatrix}-at&\tfrac{1}{4t}-a^2t\\t&at\end{bmatrix}
     \quad
 \td{A}_2 =
 \begin{bmatrix} 1/2 & a/2 \\0&0\end{bmatrix}
 \quad\text{where } t>0 \And a\in\bbC .
\]
Since this is a contractive pair,
$Z:=I-\td{A}_1\td{A}_1^*-\td{A}_2\td{A}_2^*\ge0$.  It suffices to
check that $Z_{22} \ge 0$ and $\det Z \ge 0$.  In other words,
\smallbreak
{\renewcommand{\labelenumi}{(\arabic{enumi})}
\begin{enumerate}
 \item \qquad $(1+|a|^2)t^2\le1$ \quad and
 \item \qquad $12(1+|a|^2)^2t^4-(13-4|a|^2+8\re a^2)t^2+1\le0$.
\end{enumerate}    }
\smallbreak\noindent
The dilation of this pair has pure rank 1 if and only if the
determinant is 0, which requires an equality in (2).
Notice that when $a=0$, one obtains the inequality
$1/\sqrt{12} \le t \le 1$.  The extremes yield the two pure index 1
pairs $(B_1,B_2)$ and $(C_1,C_2)$ obtained above.
%
\end{eg}


%

\newpage

\begin{tabbing}
{\it E-mail addresses}:xx\= \kill
\noindent {\footnotesize\it Address}:
\>{\footnotesize\sc Pure Mathematics Dept.}\\
\>{\footnotesize\sc University of Waterloo}\\
\>{\footnotesize\sc Waterloo, ON\quad N2L 3G1}\\
\>{\footnotesize\sc CANADA}\\
\\
{\footnotesize\it E-mail addresses}:
\>{\footnotesize\sf krdavidson@math.uwaterloo.ca}\\
\>{\footnotesize\sf dwkribs@barrow.uwaterloo.ca}\\
\>{\footnotesize\sf meshpigel@math.uwaterloo.ca}
\end{tabbing}

\end{document}